\documentclass[11pt,dina4,a4paper, twoside]{article}
\usepackage{mathbbol}
\usepackage{mathrsfs}
\usepackage{amsfonts}
\usepackage{amsmath}
\usepackage{stmaryrd}
\usepackage{booktabs}
\usepackage{cite}
\usepackage{amssymb}
\usepackage{verbatim}
\usepackage[all,cmtip]{xy}
\usepackage[top=3.54cm,bottom=4.54cm,left=2.94cm,right=2.94cm ]{geometry}

\usepackage{multirow}
\usepackage{cite}
\usepackage{caption}
\usepackage{tabularx}
\usepackage{color}
\usepackage{graphicx}
\usepackage{subfigure}

\newcommand{\p}{\partial}
\newcommand{\og}{\omega}
\newcommand{\Og}{\Omega}
\newcommand{\ep}{\varepsilon}
\newcommand{\fl}[2]{\frac{#1}{#2}}
\newcommand{\dt}{\delta}
\newcommand{\nn}{\nonumber}

\newcommand{\tht}{\theta}

\newcommand{\wt}{\widetilde}
\newcommand{\Dt}{\Delta}

\newcommand{\be}{\begin{equation}}
\newcommand{\ee}{\end{equation}}
\newcommand{\ba}{\begin{array}}
\newcommand{\ea}{\end{array}}
\newcommand{\bea}{\begin{eqnarray}}
\newcommand{\eea}{\end{eqnarray}}
\newcommand{\beas}{\begin{eqnarray*}}
\newcommand{\eeas}{\end{eqnarray*}}
\newtheorem{theorem}{Theorem}[section]
\newtheorem{lemma}{Lemma}[section]
\newtheorem{remark}{Remark}[section]

\numberwithin{equation}{section}
\newcommand{\bx}{{\bf x} }

\newcommand{\vphi}{\varphi}

\title{Uniform error bounds of a finite difference method for the Klein-Gordon-Zakharov system in the subsonic limit regime
via an asymptotic consistent formulation\thanks{This work was partially supported by the Ministry
of Education of Singapore grant
R-146-000-223-112 (W. Bao) and the Natural Science Foundation
of China Grants 91430103 and U1530401 (C. Su).}}
\author{Weizhu Bao\thanks{Department of Mathematics,
National University of Singapore, Singapore 119076 ({\tt matbaowz@nus.edu.sg},
URL: http://www.math.nus.edu.sg/\~{}bao/)} \and
Chunmei Su\thanks{Beijing Computational Science Research Center, Beijing 100193,
China; and Department of Mathematics,
National University of Singapore, Singapore 119076 ({\tt sucm@csrc.ac.cn})}}

\date{}
\begin{document}

\maketitle
\begin{abstract}
 We establish uniform error bounds of a finite
 difference method for the
 Klein-Gordon-Zakharov system (KGZ) with a
 dimensionless parameter $\ep \in (0,1]$,
 which is inversely proportional to the acoustic speed. In the subsonic limit regime, i.e. $0<\ep \ll 1$, the solution propagates highly oscillatory waves in time and/or rapid outgoing initial
layers in space due to the singular perturbation in the
Zakharov equation and/or the incompatibility of the
initial data. Specifically, the solution propagates waves with $O(\ep)$-wavelength in time and
$O(1)$-wavelength in space as well as outgoing initial layers in space
at speed $O(1/\ep)$. This high oscillation  in time and rapid outgoing waves in space of the solution cause
significant burdens in designing numerical methods and establishing error estimates for
KGZ. By adapting an asymptotic consistent formulation, we propose a
uniformly accurate finite difference method and rigorously establish
two independent error bounds at $O(h^2+\tau^2/\ep)$ and
$O(h^2+\tau+\ep)$ with $h$ mesh size and $\tau$ time step. Thus we
obtain a uniform error bound at $O(h^2+\tau)$ for $0<\ep\le 1$.
The main techniques in the analysis include the energy
method, cut-off of the nonlinearity to bound the numerical solution, the integral approximation of
the oscillatory term, and $\ep$-dependent error bounds
between the solutions of KGZ and its limiting model when $\ep\to0^+$.
Finally, numerical results are
reported to confirm our error bounds.
\end{abstract}



\pagestyle{myheadings}\thispagestyle{plain}

% =============================================================================
%                                                                              SECTION 1
% =============================================================================
\section{Introduction}
\label{section1}
We study the Klein-Gordon-Zakharov (KGZ) system which describes the interaction between  a Langmuir wave and an ion acoustic wave in plasma \cite{Masmoudi}:
\be\label{KGZo}
\begin{aligned}
& \p_{tt}E(\bx,t)-3v_0^2\Dt E(\bx,t)+\og_p^2 E(\bx,t)+\og_p^2 N(\bx,t)E(\bx,t)=0,\quad \bx\in\mathbb{R}^d,\quad t>0,\\
&\p_{tt} N(\bx,t)-c_s^2 \Dt N(\bx,t)-\frac{n_0\ep_0}{4mN_0}\Dt |E|^2(\bx,t)=0, \quad \bx \in \mathbb{R}^d, \quad  t>0,
%&E^\ep(\bx,0)=E_0(\bx),\quad \p_t E^\ep(\bx,0)=E_1(\bx),\quad \bx\in\mathbb{R}^d,\\
%&N^\ep(\bx,0)=N_0^\ep(\bx),\quad \p_t N^\ep (\bx,0)=N_1^\ep(\bx),\quad \bx\in\mathbb{R}^d.
\end{aligned}
\ee
where $t$ is time, $\bx\in {\mathbb R}^d$ ($d=1,2,3$) is the spatial coordinate, $E(\bx,t)$ and $N(\bx,t)$ are real-valued functions representing the fast time scale component of the electric field raised by electrons and the ion density fluctuation from the constant equilibrium, respectively. Here $v_0$ is the electron thermal velocity, $\og_p$ is the electron plasma frequency, $c_s$ is the acoustic speed, $n_0$ is plasma charge number, $\ep_0$ is vacuum dielectric constant, $m$ is ion mass and $N_0$ is electron density. It can be derived from the Euler equations for the electrons and ions, coupled with the Maxwell equation for the electron field by negeleting the magnetic effect and further assuming that ions move much slower than electrons (cf. \cite{Berge, Dendy, Sulem, Zakharov} for physical and formal derivations and \cite{Texier} for mathematical justifications).

For scaling the KGZ system \eqref{KGZ}, we introduce
\be\label{dimen}
\wt{t}=\frac{t}{t_s},\quad \wt{\bx}=\frac{\bx}{x_s},\quad \wt{E}(\wt{\bx},\wt{t})=\frac{E(\bx,t)}{E_s},\quad \wt{N}(\wt{\bx},\wt{t})=\frac{N(\bx,t)}{N_s},
\ee
where $t_s=\fl{1}{\og_p}$, $x_s=\fl{\sqrt{3}v_0}{\og_p}$, $E_s=2c_s \sqrt{\frac{mN_0}{n_0\ep_0}}$ and $N_s=1$ are the dimensionless time, length, electric field and ion density unit, respectively.
Plugging \eqref{dimen} into \eqref{KGZ} and removing all $\,\wt{\,\,}\,$ followed by replacing $N(\bx,t)$ and $E(\bx,t)$ by $N^\ep(\bx,t)$ and $E^\ep(\bx,t)$, respectively, we get the following dimensionless KGZ system as
\be\label{KGZ}
\begin{aligned}
&\p_{tt}E^\ep(\bx,t)-\Dt E^\ep(\bx,t)+E^\ep(\bx,t)+N^\ep(\bx,t)E^\ep(\bx,t)=0,\quad \bx\in\mathbb{R}^d,\quad t>0,\\
&\ep^2\p_{tt} N^\ep(\bx,t)-\Dt N^\ep(\bx,t)-\Dt |E^\ep|^2(\bx,t)=0, \quad \bx \in \mathbb{R}^d, \quad  t>0,
\end{aligned}
\ee
where the dimensionless parameter $0<\ep:=\fl{\sqrt{3}v_0}{c_s}\le 1$  is inversely proportional to the speed of sound.
%\be\label{epg}
%\ep.
%\ee
Here we consider the case where the thermal electron velocity is much smaller than the ion-acoustic speed, i.e. $3v_0^2\ll c_s^2$, which gives $0<\ep\ll 1$, i.e. the KGZ system in the subsonic limit regime. To study the dynamics of the KGZ system \eqref{KGZ}, the initial data is usually given as
\be\label{initiald}
E^\ep(\bx,0)=E_0(\bx),\,\,\,\, \p_t E^\ep(\bx,0)=E_1(\bx),\,\,\, \, N^\ep(\bx,0)=N_0^\ep(\bx),\,\,\,\,
\p_t N^\ep(\bx,0)=N_1^\ep(\bx).
\ee

As is known, \eqref{KGZ} is time symmetric or time reversible and conserves the
total {\sl energy} \cite{Masmoudi, Masmoudi2008}, i.e. for $t\ge 0$
$$\mathcal {H}^\ep(t):=\int_{\mathbb{R}^d}\left[|\p_tE^\ep|^2+|\nabla E^\ep|^2+|E|^2+\fl{\ep^2}{2}
|\nabla \varphi^\ep|^2+\fl{1}{2}|N^\ep|^2+N^\ep|E^\ep|^2\right]d\bx\equiv \mathcal {H}^\ep(0),$$
where $\varphi^\ep$ is defined by $\Dt\varphi^\ep=\p_t N^\ep$ with $\lim\limits_{|\bx|\rightarrow \infty}\vphi^\ep=0$.

There have been extensive studies for the KGZ system in the
literatures for $\ep=\ep_0$ with $\ep_0>0$ a fixed constant, i.e. $O(1)$-acoustic-speed regime. Along the analytical
part, for the derivation of the KGZ from two-fluid Euler-Maxwell
system, we refer to \cite{Berge, Texier}; and for the well-posedness of the Cauchy problem, we refer to \cite{Ozawa95, Ozawa99, Tsutaya, Ohta}.
Along the numerical part, we refer to \cite{Wang} for finite difference method and \cite{Zhao2013, Zhao2015} for exponential-wave-integrator Fourier pseudospectral method. However, in the subsonic limit regime,
the analysis and efficient computation of the KGZ system are rather
complicated \cite{Berge, Masmoudi} due to the high oscillation in time
and/or rapid outgoing waves in space of the solution as $\ep\rightarrow 0^+$.

\begin{figure}[t!]
\begin{minipage}[t]{15cm}
\centering
\includegraphics[width=15cm,height=5cm]{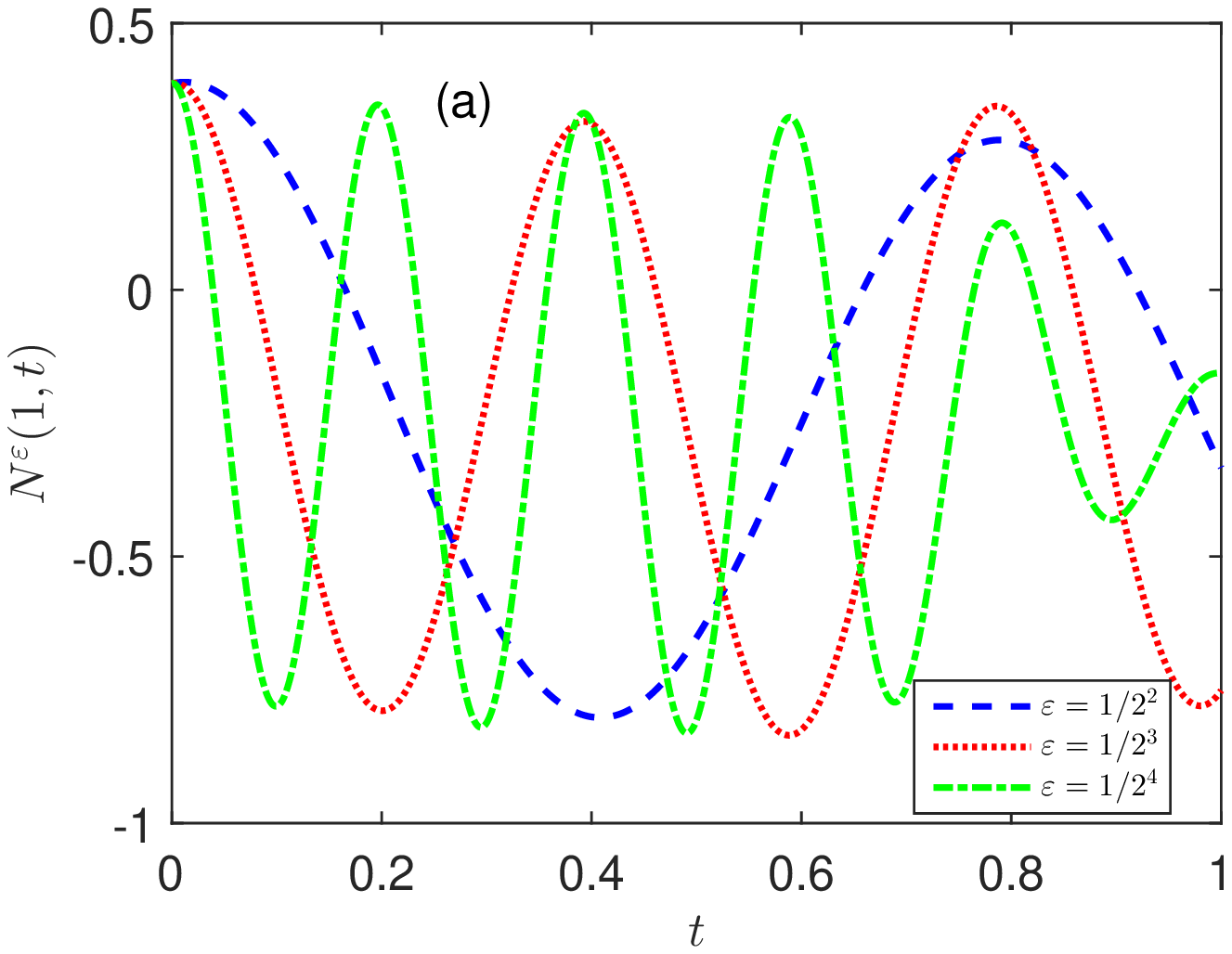}
\end{minipage}
\begin{minipage}[t]{15cm}
\centering
\includegraphics[width=15cm,height=5cm]{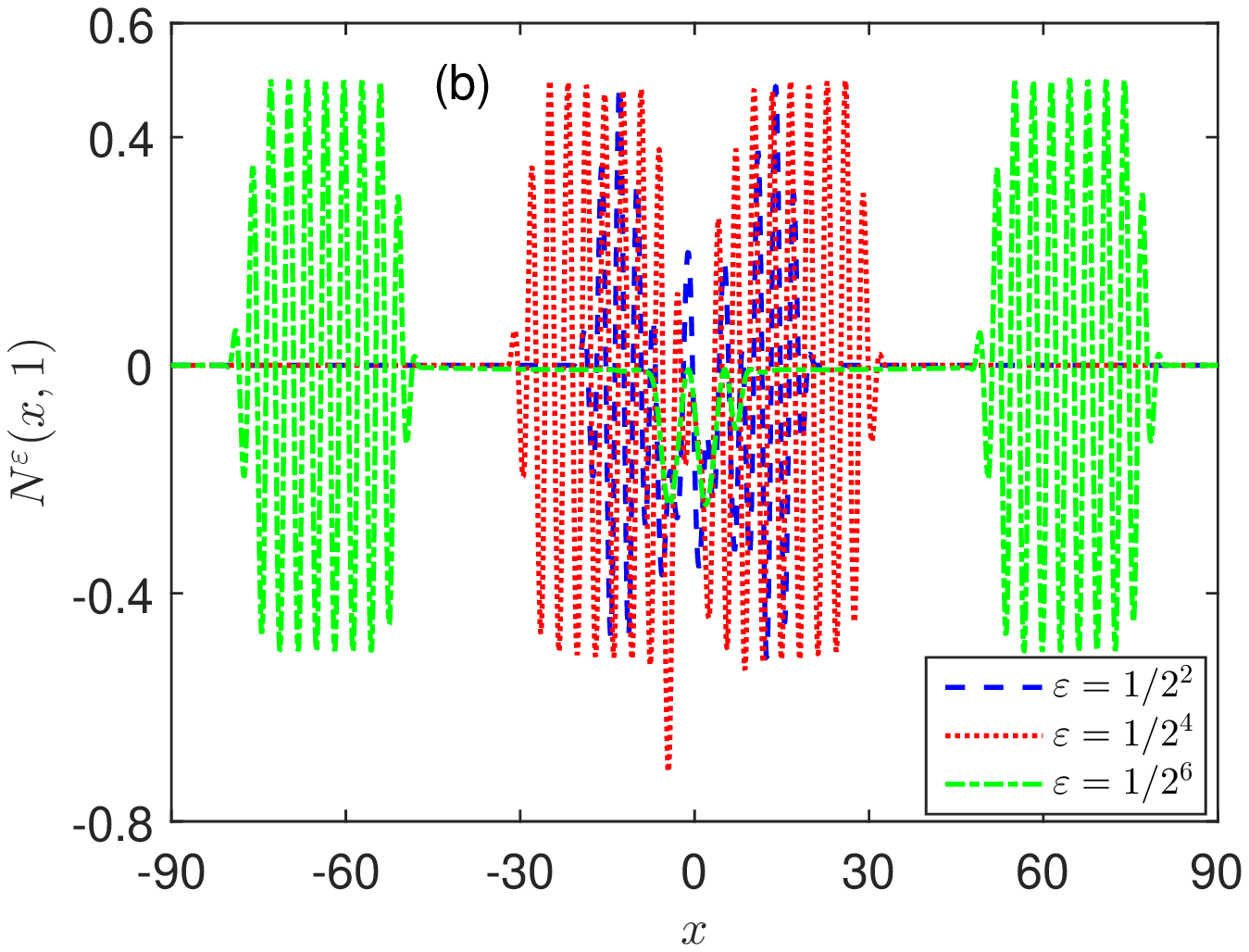}
\end{minipage}
\caption{The temporal oscillation (a) and rapid outgoing wave in space (b) of the KGZ \eqref{KGZ} for $d=1$.}\label{fig1}
\end{figure}

Based on the results in \cite{Masmoudi2010, Daub}, in the subsonic
limit, i.e. $\ep \rightarrow 0^+$, the KGZ system collapses to the
Klein-Gordon (KG) equation. Formally we have $E^\ep\rightarrow E_{\rm k}$,
where $E_{\rm k}:=E_{\rm k}(\bx,t)$ is the solution of the KG equation \cite{Masmoudi2010, Daub}:
 \be\label{KG}
\begin{aligned}
 &\p_{tt}E_{\rm k}(\bx,t)-\Dt E_{\rm k}(\bx,t)+E_{\rm k}(\bx,t)-E_{\rm k}(\bx,t)^3=0,\quad \bx\in \mathbb{R}^d, \quad t>0,\\
 & E(\bx,0)=E_0(\bx),\quad \p_t E(\bx,0)=E_1(\bx),\quad \bx \in \mathbb{R}^d.\\
  \end{aligned}
 \ee
The KG \eqref{KG} conserves the {\sl energy}
\[\mathcal{H}(t):=\int_{\mathcal{R}^d}\left[|\p_t E_{\rm k}|^2+|\nabla E_{\rm k}|^2+|E_{\rm k}|^2-\fl 12|E_{\rm k}|^4\right] d\bx\equiv \mathcal{H}(0),\quad t\ge 0.\]
Different convergence rates can be obtained due to the incompatibility of the initial data ($E_0,E_1,N_0^\ep,N_1^\ep$) for \eqref{KGZ} with respect to \eqref{KG}, which can be characterized as
\be\label{inct}
N_0^\ep(\bx)=-E_0(\bx)^2+\ep^\alpha \og_0(\bx),\quad N_1^\ep(\bx)=-2E_0(\bx)E_1(\bx)+\ep^\beta \og_1(\bx),\quad \bx\in\mathbb{R}^d,
\ee
where $\alpha\ge 0$ and $\beta\ge -1$ are given parameters and $\og_0(\bx)$ and  $\og_1(\bx)$ are
given functions, which are all independent of $\ep$. Similar to the properties of the solutions of the Zakharov system in the subsonic limit regime \cite{Masmoudi2008, Ozawa,Schochet}, when $0<\ep\ll 1$,
the solution of the KGZ system propagates waves with wavelength $O(\ep)$
and $O(1)$ in time and space, respectively (cf. Fig. \ref{fig1}a), and/or rapid outgoing initial layers at speed $O(1/\ep)$ in space (cf. Fig. \ref{fig1}b). More precisely, when $\alpha\ge 2$ and $\beta\ge 1$, the
leading oscillation comes from the $\ep^2\p_{tt}$ term; and otherwise from the incompatibility of the initial data.
To illustrate the temporal oscillation and rapid outgoing wave phenomena,
Figure \ref{fig1} shows the solutions $N^\ep(x,1)$, $N^\ep(1,t)$ of the
KGZ \eqref{KGZ} for $d=1$ and the initial data
\begin{align}
&E_0(x)=\fl{1}{2}\psi\left(\fl{x+15}{8}\right)\psi\left(\fl{15-x}{7}\right)\cos \left(\fl{x}{2}\right),\ E_1(x)=\fl{1}{2}\psi\left(\fl{x+10}{5}\right)\psi\left(\fl{10-x}{5}\right)\sin \left(\fl{x}{2}\right),\nonumber\\
&\og_0(x)=\psi\left(\fl{x+18}{10}\right)\psi\left(\fl{18-x}{9}\right)\sin \left(2x+\fl{\pi}{6}\right),\quad \og_1(x)=e^{-x^2/3}\sin(2x),\label{inite}
\end{align}
with
\be\label{fg}
\psi(x)=\fl{\varphi(x)}{\varphi(x)+\varphi(1-x)}, \quad \varphi(x)=e^{-1/x}\chi_{(0,\infty)},
\ee
and $\chi$ being the characteristic function, $\alpha=\beta=0$ in \eqref{inct} for different $\ep$,
which was obtained numerically by the exponential-wave-integrator sine pseudospectral method on a bounded interval
$[-200,200]$ with the homogenous Dirichlet boundary condition \cite{Zhao2013}.

The highly temporal oscillatory nature in the solution of the KGZ \eqref{KGZ} brings significant numerical difficulties, especially in the subsonic limit regime, i.e. $0<\ep\ll 1$. To the best of our knowledge, there are few results concerning error estimates of different numerical methods for KGZ with respect to mesh size $h$ and time step $\tau$ as well as the parameter $0<\ep\le 1$. Recently, a conservative finite difference method (FDM) was proposed and analyzed in the subsonic limit regime \cite{KGZ1}, where it was proved that in order to obtain `correct' oscillatory solutions, the FDM requests the meshing strategy (or $\ep$-scalability) $h=O(\ep^{1/2})$ and $\tau=O(\ep^{3/2})$. The reason is due to that
$N^\ep(\bx,t)$ does not converge as $\ep\to0^+$ when $\alpha=0$ or $\beta=-1$ \cite{Masmoudi2008,Schochet,Sulemb} (cf. Figure \ref{fig1}).

The main aim of this paper is to propose and analyze a finite difference method for KGZ,
which is uniformly accurate in both space and time for $0<\ep\le1$. The key points in
designing the uniformly accurate finite difference method include (i) reformulating KGZ
into an asymptotic consistent formulation and (ii) adopting an integral approximation of
the oscillatory term. To establish the error bounds, we apply the energy method, cut-off
technique for treating the nonlinearity and the inverse estimates to bound the numerical
solution, and the limiting equation via a nonlinear Klein-Gordon equation with an oscillatory potential. The error bounds of our new numerical method significantly relax
the meshing strategy of the standard FDM for KGZ in the subsonic
limit regime \cite{KGZ1}.

The rest of the paper is organized as follows. In section 2, we introduce an asymptotic
consistent formulation of KGZ, present a finite difference method and state our main results. Section 3 is devoted to the details of the error analysis. Numerical results are reported
in section 4 to confirm our error bounds. Finally some conclusions are drawn in section 5.
Throughout the paper, we adopt the standard Sobolev spaces as well as the corresponding
norms and denote $A\lesssim B$ to represent that there exists a generic constant $C>0$ independent of $\ep$, $\tau$, $h$, such that $|A|\le C\,B$.

\section{A finite difference method and its error bounds}
In this section, we present a uniformly accurate finite difference method based on an asymptotic consistent formulation of KGZ and give its uniform error bounds.
\subsection{An asymptotic consistent formulation}
Following \cite{Bao2016}, we introduce
\be\label{Nd}
F^\ep(\bx,t)=N^\ep(\bx,t)+|E^\ep(\bx,t)|^2-G^\ep(\bx,t), \quad \bx\in \mathbb{R}^d,\quad t\ge0,
\ee
where $G^\ep(\bx,t)$ represents initial layer caused by the incompatibility of the initial data \eqref{inct}, which is the solution of the linear wave equation
\be\label{wave}
\begin{split}
&\p_{tt}G^\ep(\bx,t)-\fl{1}{\ep^2}\Delta G^\ep(\bx,t)=0, \quad \bx\in \mathbb{R}^d,\quad t>0,\\
&G^\ep(\bx,0)=\ep^\alpha \og_0(\bx), \quad \p_t G^\ep(\bx,0)=\ep^\beta \og_1(\bx).
\end{split}
\ee
Substituting \eqref{Nd} into the KGZ \eqref{KGZ}, we can reformulate it into an asymptotic consistent formulation
\be\label{KGZ2}
     \begin{aligned}
     & \p_{tt}E^\ep(\bx,t)-\Dt E^\ep(\bx,t)+\left[1-E^\ep(\bx,t)^2+F^\ep(\bx,t)+ G^\ep(\bx,t)\right]E^\ep(\bx,t)=0,  \\
     &\p_{tt} F^\ep(\bx,t)-\frac{1}{\ep^2}\Dt F^\ep(\bx,t)-\p_{tt} |E^\ep(\bx,t)|^2=0, \quad \bx \in \mathbb{R}^d, \quad t>0,\\
     &E^\ep(\bx,0)=E_0(\bx),\quad \p_t E^\ep(\bx,0)=E_1(\bx),\quad F^\ep(\bx,0)=0,\quad \p_t F^\ep (\bx,0)=0.
                \end{aligned}
\ee

\begin{figure}[t!]
\begin{minipage}[t]{5cm}
\centering
\includegraphics[width=5cm,height=5.2cm]{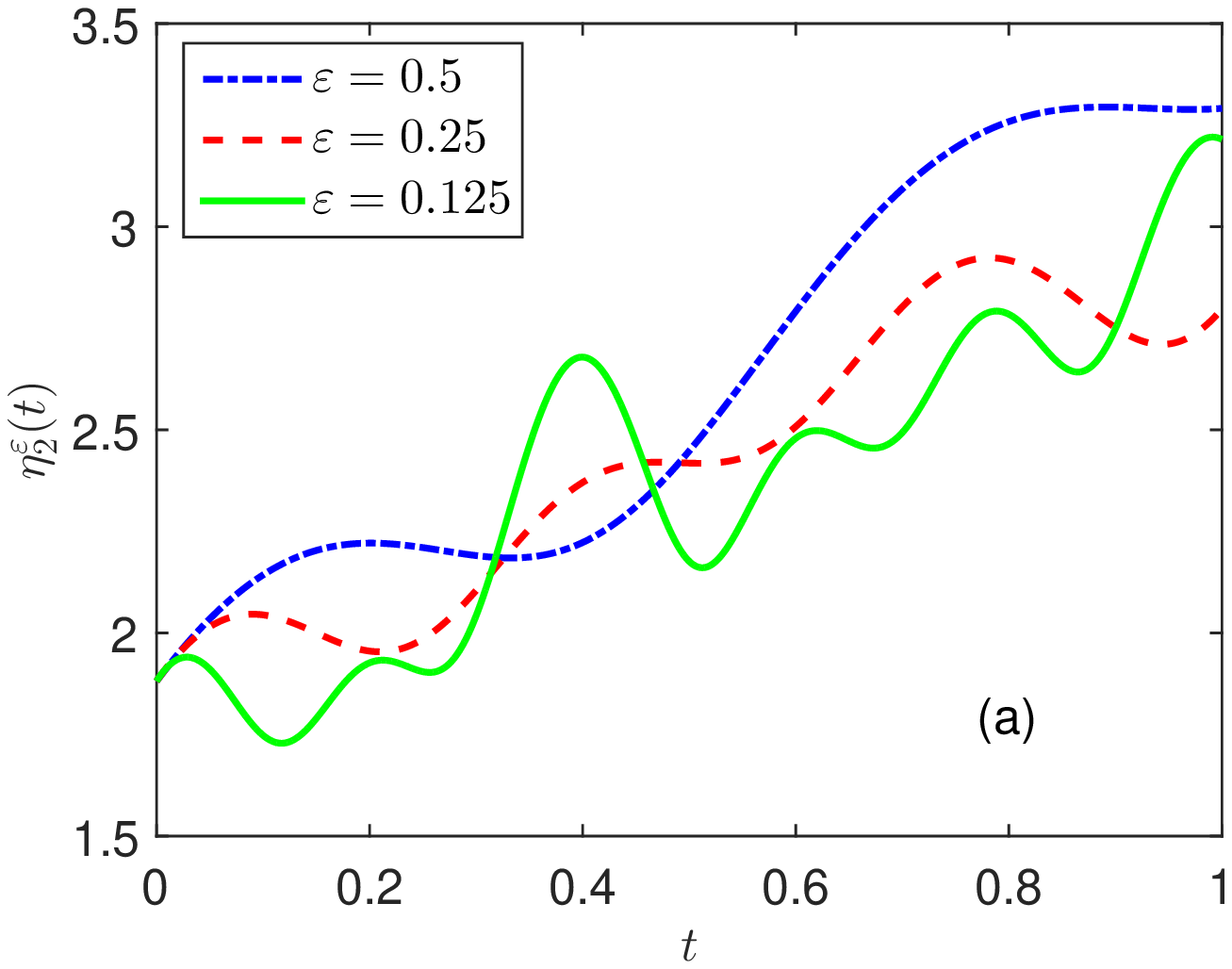}
\end{minipage}
\hspace{-3mm}
\begin{minipage}[t]{5cm}
\centering
\includegraphics[width=5cm,height=5.2cm]{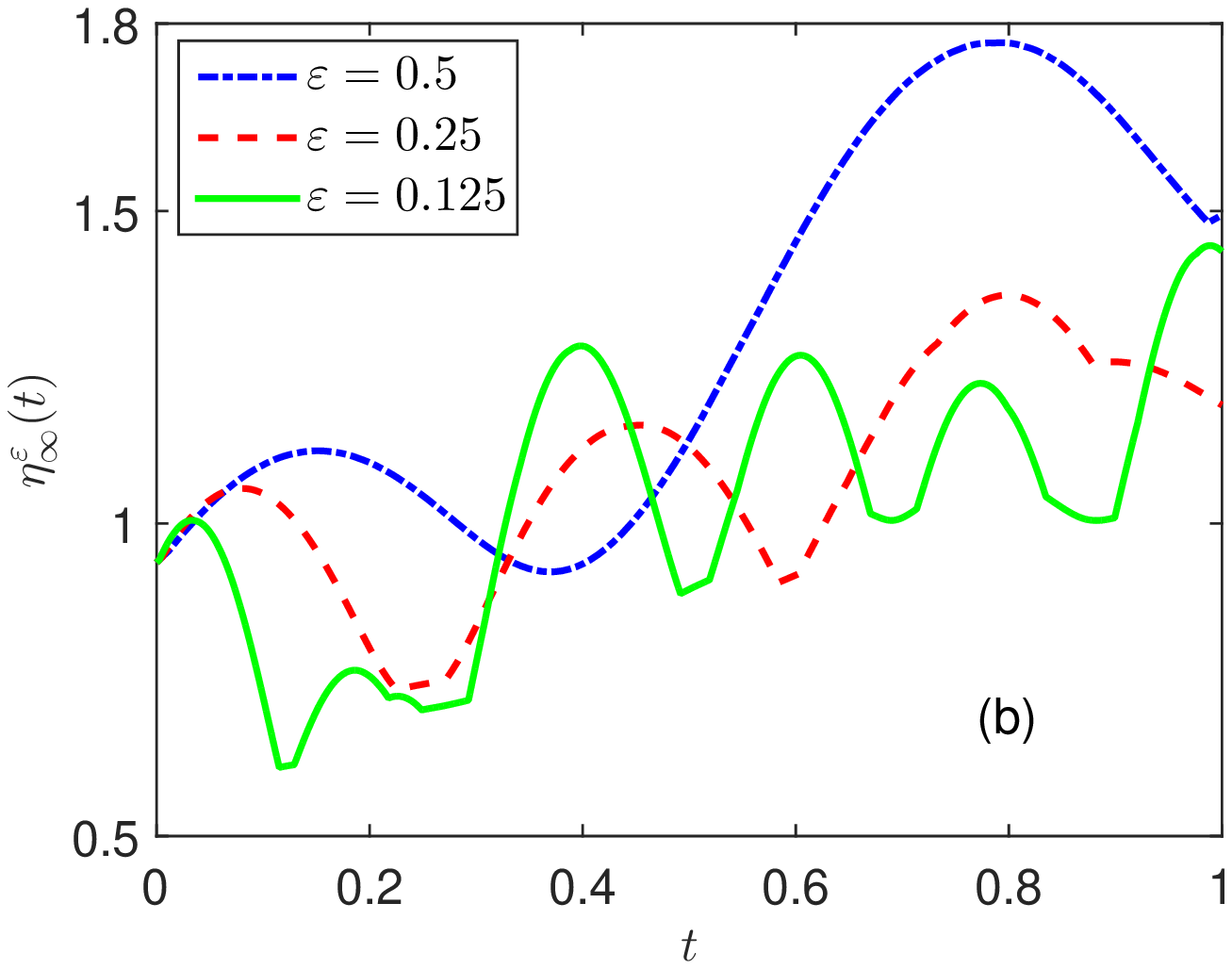}
\end{minipage}
\hspace{-3mm}
\begin{minipage}[t]{5cm}
\centering
\includegraphics[width=5cm,height=5.2cm]{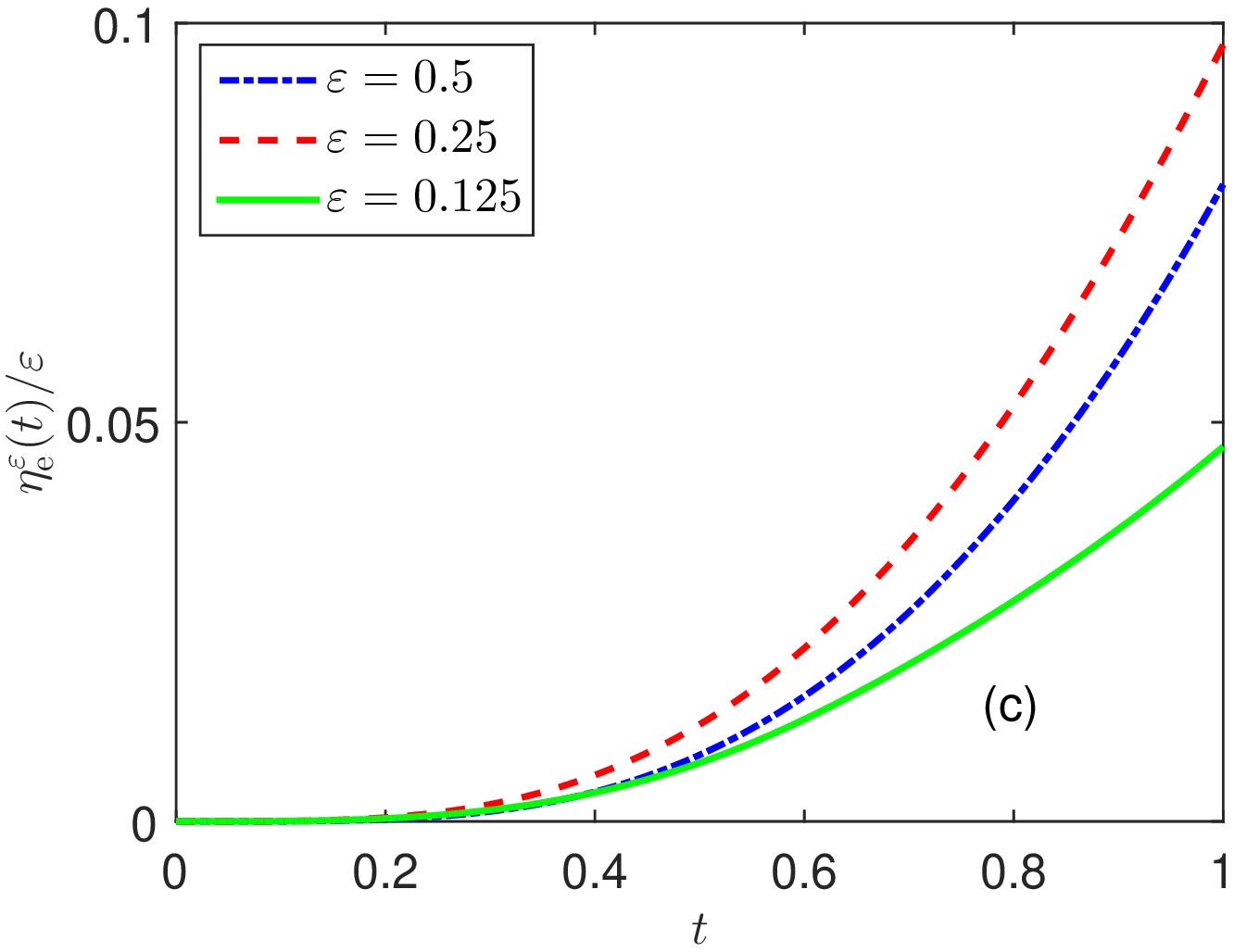}
\end{minipage}
\caption{Time evolution of $\eta_2^\ep(t)$, $\eta_\infty^\ep(t)$ and
$\eta_{\rm e}^\ep(t)$.}\label{fig2}
\end{figure}

In the subsonic limit regime, i.e. $\ep \to 0^+$, formally
we have $E^\ep(\bx,t)\to E_{\rm k}(\bx,t)$ and $F^\ep(\bx,t)\to 0$, where $E_{\rm k}(\bx,t)$
is the solution of the KG \eqref{KG}. Moreover, as $\ep\to0^+$, formally we can also get
$E^\ep(\bx,t)\to \wt{E}^\ep(\bx,t)$, where $\wt{E}^\ep:=\wt{E}^\ep(\bx,t)$ is the solution of the Klein-Gordon equation with an oscillatory potential $G^\ep(\bx,t)$ (KG-OP)
\be\label{AKG}
     \begin{aligned}
     & \p_{tt}\wt{E}^\ep(\bx,t)-\Delta \wt{E}^\ep(\bx,t)+\left[1-\wt{E}^\ep(\bx,t)^2+ G^\ep(\bx,t)\right]\wt{E}^\ep(\bx,t)=0, \\ &\wt{E}^\ep(\bx,0)=E_0(\bx),\quad
      \p_t\wt{E}^\ep(\bx,0)=E_1(\bx), \quad\bx\in\mathbb{R}^d.
                \end{aligned}
\ee
Inspired by the convergence of the Zakharov system to the nonlinear Schr\"{o}dinger equation in the subsonic limit regime \cite{Ozawa} and the analytical analysis of the KGZ converging to the KG \cite{Daub}, we can obtain the following result concerning on the convergence from the KGZ \eqref{KGZ} to the KG-OP \eqref{AKG}
\be\label{conv}
\|F^\ep\|_{L^2}+\|F^\ep\|_{L^\infty}+\|E^\ep(\cdot,t)-\wt{E}^\ep(\cdot,t)\|_{H^1}\le C_T \ep, \quad 0\le t\le T,
\ee
where $0<T<T^*$ with $T^*>0$ being the maximum common existence time for the solutions of the KGZ \eqref{KGZ} and the KG-OP \eqref{AKG} and $C_T$ is a positive constant independent of $\ep$. To illustrate this, Figure \ref{fig2} depicts the convergence behavior between the solutions of the KGZ \eqref{KGZ} and the KG-OP \eqref{AKG}, where
$\eta_2^\ep(t):=\fl{1}{\ep}\|F^\ep(\cdot,t)\|_{L^2}+\|\p_t F^\ep(\cdot,t)\|_{L^2}+\|\p_{tt} F^\ep(\cdot,t)\|_{L^2}$,
$\eta_\infty^\ep(t):=\fl{1}{\ep}\|F^\ep(\cdot,t)\|_{L^\infty}+\|\p_t F^\ep(\cdot,t)\|_{L^\infty}+\|\p_{tt} F^\ep(\cdot,t)\|_{L^\infty}$ and
$\eta_{\rm e}^\ep(t):=\|E^\ep(\cdot,t)-\wt{E}^\ep(\cdot,t)\|_{H^1}$
for different $\ep$ with the same initial data as in \eqref{inite} for $d=1$ and $\alpha=\beta=0$.

\subsection{A uniformly accurate finite difference method}
For simplicity of notation, we will only present the numerical method for the KGZ system
on one space dimension, and the extensions to higher dimensions are straightforward.
Practically, similar to most works for computation of the Zakharov-type equations \cite{Bao2005,Jin,Cai2015,Bao2016,KGZ1}, \eqref{KGZ2} is
truncated on a bounded interval $\Og=(a,b)$ with the homogeneous Dirichlet boundary condition:
\be\label{KGZ1d}
        \begin{split}
    &\p_{tt}E^\ep(x,t)-\p_{xx} E^\ep(x,t)+\left[1-E^\ep(x,t)^2+F^\ep(x,t)+ G^\ep(x,t)\right]E^\ep(x,t)=0,  \\
     &\p_{tt} F^\ep(x,t)-\frac{1}{\ep^2}\p_{xx} F^\ep(x,t)-\p_{tt} |E^\ep(x,t)|^2=0,  \quad x \in \Omega, \quad t>0,\\
     &E^\ep(x,0)=E_0(x),\,\,\,\, \p_t E^\ep(x,0)=E_1(x),\,\,\,\, F^\ep(x,0)=0,\,\,\,\, \p_t F^\ep (x,0)=0, \,\,\,\, x\in \overline\Omega,\\
     &E^\ep(a,t)=E^\ep(b,t)=0,\quad F^\ep(a,t)=F^\ep(b,t)=0, \quad t\ge0,
                \end{split}
\ee
where $G^\ep(x,t)$ is defined as \eqref{wave} with homogeneous Dirichlet boundary condition for $d=1$,
\be\label{wave1D}
\begin{split}
&\p_{tt}G^\ep(x,t)-\fl{1}{\ep^2}\partial_{xx} G^\ep(x,t)=0, \quad x\in \Omega,\quad t>0,\\
&G^\ep(x,0)=\ep^\alpha \og_0(x), \,\, \p_t G^\ep(x,0)=\ep^\beta \og_1(x), \,\,
x\in \overline{\Omega};\,\,G^\ep(a,t)=G^\ep(b,t)=0, \,\, t\ge0.
\end{split}
\ee

As $\ep\to 0$, formally we have $E^\ep(x,t)\to \wt{E}^\ep(x,t)$ and $F^\ep(x,t)\to 0$, where $\wt{E}^\ep(x,t)$ is the solution of the KG-OP
\be\label{AKG1d}
         \begin{aligned}
     & \p_{tt}\wt{E}^\ep(x,t)-\p_{xx} \wt{E}^\ep(x,t)+\left[1-\wt{E}^\ep(x,t)^2+ G^\ep(x,t)\right]\wt{E}^\ep(x,t)=0,  \\
     &\wt{E}^\ep(x,0)=E_0(x),\,\,\,\, \p_t \wt{E}^\ep(x,0)=E_1(x),\,\,\,\, x\in \overline{\Omega};\,\,\,\,\wt{E}^\ep(a,t)=\wt{E}^\ep(b,t)=0, \,\,\,\, t\ge0.
                \end{aligned}
      \ee

Choose a mesh size $h:=\Dt x=\fl{b-a}{M}$ with $M$ being a positive integer and a time step $\tau:=\Dt t>0$ and denote the grid points and time steps as
$$x_j:=a+jh,\quad j=0,1,\cdots,M;\quad t_k:=k\tau,\quad k=0,1,2,\cdots.$$
Define the index sets
$$\mathcal {T}_M=\{j \ | \ j=1,2,\cdots,M-1\},\quad \mathcal{T}_M^0=\{j\ |\ j=0,1,\cdots,M\}.$$
Let $E^{\ep,k}_j$ and $F^{\ep,k}_j$ be the approximations of $E^{\ep}(x_j,t_k)$ and $F^\ep(x_j,t_k)$, respectively, and denote $E^{\ep,k}=(E_0^{\ep,k}, E_1^{\ep,k},\cdots, E_M^{\ep,k})^T$, $F^{\ep,k}=(F_0^{\ep,k}, F_1^{\ep,k},\cdots, F_M^{\ep,k})^T\in \mathbb{R}^{(M+1)}$ as the numerical solution vectors at $t=t_k$. The finite difference operators are the standard notations as:
\beas
&&\dt_x^+E_j^k=\fl{E_{j+1}^k-E_j^k}{h},\quad
\dt_t^+E_j^k=\fl{E_j^{k+1}-E_j^k}{\tau},\quad
\dt_t^c E_j^k=\fl{E_j^{k+1}-E_j^{k-1}}{2\tau},\\
&&\dt_t^2E_j^k=\fl{E_j^{k+1}-2E_j^k+E_j^{k-1}}{\tau^2},\quad
\dt_x^2E_j^k=\fl{E_{j+1}^k-2E_j^k+E_{j-1}^k}{h^2}.
\eeas

In this paper, we consider the finite difference discretization of \eqref{KGZ1d} as following
\be\label{scheme1}
\begin{aligned}
&\dt_t^2E_j^{\ep,k}=\left(\dt_x^2-1+|E_j^{\ep,k}|^2-F_j^{\ep,k}-
H_j^{\ep,k}\right)\fl{E_j^{\ep,k+1}+E_j^{\ep,k-1}}{2},\\
&\dt_t^2F_j^{\ep,k}=\fl{1}{2\ep^2}\dt_x^2(F_j^{\ep,k+1}+F_j^{\ep,k-1})+
\dt_t^2|E_j^{\ep,k}|^2, \quad j\in \mathcal{T}_M, \quad k\ge 1.
\end{aligned}
\ee
where we apply an average of the oscillatory potential $G^\ep$ over the interval $[t_{k-1},t_{k+1}]$
\be\label{Hd}
H_j^{\ep,k}=\int_{-1}^{1}(1-|s|)G^\ep(x_j,t_k+s\tau)ds,\quad j\in \mathcal{T}_M,\quad k\ge 1.
\ee
 Meanwhile, the boundary and initial conditions are discretized as
\be\label{initfd11}
E_0^{\ep,k}=E_M^{\ep,k}=F_0^{\ep,k}=F_M^{\ep,k}=0,\quad k\ge0; \quad
E_j^{\ep,0}=E_0(x_j),\quad F_j^{\ep,0}=0, \quad j\in \mathcal{T}_M^0.
\ee
Next we consider the value of the first step $E_j^{\ep,1}$ and $F_j^{\ep,1}$. By Taylor expansion, we get $E_j^{\ep,1}$ as
\be\label{E_1}
E_j^{\ep,1}=E_0(x_j)+\tau E_1(x_j)+\fl{\tau^2}{2}\p_{tt}
E^\ep(x_j,0),\quad F_j^{\ep,1}=\fl{\tau^2}{2}\p_{tt}F^\ep(x_j,0),\quad j\in \mathcal{T}_M,
\ee
where by \eqref{KGZ1d},
\[
\p_{tt}E^\ep(x,0)=E_0''(x)-E_0(x)-N_0^\ep(x)E_0(x),\quad
\p_{tt}F^\ep(x,0)=2E_1(x)^2+2E_0(x)\p_{tt}E^\ep(x,0).
\]

In practical computation, $H_j^{\ep,k}$ in \eqref{Hd} can be obtained by solving the wave equation \eqref{wave1D} via the sine pseudospectral discretization in space followed by integrating in time in phase space {\sl exactly} \cite{Bao2016} as
\begin{align*}
H_j^{\ep,k}&=\ep^\alpha\sum_{l=1}^{M-1}\wt{(\omega_0)}_l \sin
\left(\frac{jl\pi}{M}\right)\int_{-1}^1 (1-|s|)
\cos\left(\tht_l (t_k+s\tau)\right)ds \\
&\quad +\ep^\beta \sum_{l=1}^{M-1}\fl{\wt{(\omega_1)}_l}{\tht_l} \sin
\left(\frac{jl\pi}{M}\right)\int_{-1}^1 (1-|s|)
\sin\left(\tht_l (t_k+s\tau)\right)ds \\
&=2\sum_{l=1}^{M-1} \sin
\left(\frac{jl\pi}{M}\right)
\left[\ep^\alpha \wt{(\omega_0)}_l\cos\left(\tht_l t_k\right)+\ep^{\beta}\frac{\wt{(\omega_1)}_l}{\tht_l}
\sin\left(\tht_l t_k\right)\right]\int_0^1
\cos\left(\tht_l s\tau\right)(1-s)ds,\\
&=\fl{4}{\tau^2}\sum_{l=1}^{M-1}\frac{1}{\tht_l^2}\sin
\left(\frac{jl\pi}{M}\right)\sin^2\left(\frac{\tht_l\tau}{2}\right)
\left[\ep^\alpha \wt{(\omega_0)}_l\cos\left(\tht_l t_k\right)+\ep^{\beta}\frac{\wt{(\omega_1)}_l}{\tht_l}
\sin\left(\tht_l t_k\right)\right],
\end{align*}
where for $l\in \mathcal{T}_M$,
\[
\tht_l=\frac{l\pi}{\ep(b-a)}, \quad \wt{(\omega_0)}_l=\frac{ 2}{M}\sum_{j=1}^{M-1}\omega_0(x_j)
\sin\left(\frac{jl\pi}{M}\right), \quad \wt{(\omega_1)}_l=\frac{2}{M}\sum_{j=1}^{M-1}\omega_1(x_j)
\sin\left(\frac{jl\pi}{M}\right).
\]
\subsection{Main results}
For simplicity of notation, we denote
\[\alpha^*:=\min\{\alpha,1+\beta\}\ge 0.\]
Let $T^*>0$ be the maximum common existence time for the solution to the KGZ \eqref{KGZ1d} and the KG-OP \eqref{AKG1d}. Then
for $0<T<T^*$, according to the known results in \cite{Added, Masmoudi2008,Ozawa, Schochet}, we can assume the exact solution
$(E^\ep(x,t), F^\ep(x,t))$ of the KGZ \eqref{KGZ1d} and the solution
$\wt{E}^\ep(x,t)$ of the KG-OP \eqref{AKG1d} are smooth enough and satisfy
\begin{equation*}
(\rm A)\begin{split}
&\|E^\ep\|_{W^{4,\infty}(\Omega)}+\|\p_t E^\ep\|_{W^{2,\infty}(\Omega)}+
\|\p_{tt}E^\ep\|_{W^{2,\infty}(\Omega)}+\ep\|\p_t^3 E^\ep\|_{W^{2,\infty}(\Omega)}\lesssim1,\\
&\|\wt{E}^\ep\|_{W^{4,\infty}(\Omega)}+\|\p_t\wt{E}^\ep\|_{W^{2,\infty}(\Omega)}+
\|\p_{tt}\wt{E}^\ep\|_{W^{2,\infty}(\Omega)}\lesssim 1,\quad  \|\p_t^3\wt{E}^{\ep}\|_{L^\infty(\Omega)}\lesssim \fl{1}{\ep^{1-\alpha^*}},\\
&\|F^\ep\|_{W^{4,\infty}(\Omega)}\lesssim \ep,\quad \|\p_t F^\ep\|_{W^{4,\infty}(\Omega)}
+\|\p_{tt}F^\ep\|_{W^{2,\infty}(\Omega)}+\ep\|\p_t^3 F^\ep\|_{W^{2,\infty}(\Omega)}\lesssim 1.
\end{split}
\end{equation*}
Furthermore, we assume that the initial data satisfies
\begin{equation*}
(\rm B)\hskip3cm
\|E_0\|_{W^{5,\infty}(\Omega)}+\|E_1\|_{W^{5,\infty}(\Omega)}+\|\og_0\|_{W^{3,\infty}(\Omega)}+
\|\og_1\|_{W^{3,\infty}(\Omega)}\lesssim 1.\hskip7cm
\end{equation*}
It can be concluded from \eqref{wave} and assumption (B) that
\be\label{Gp}
\|\p_t^m G^\ep\|_{W^{3,\infty}(\Omega)}\lesssim \ep^{\alpha^*-m},\quad m=0,1,2,3.
\ee

To measure the error between the exact solution and the numerical solution of the KGZ system, we introduce some notations. Denote
\[X_M=\{v=(v_j)_{j\in\mathcal{T}_M^0} | \,v_0=v_M=0\} \subseteq \mathbb{R}^{M+1}.\]
The norms and inner products over $X_M$ are defined as
\beas
&&\|u\|^2=h\sum\limits_{j=1}^{M-1}|u_j|^2, \quad\|\dt_x^+ u\|^2=h\sum\limits_{j=0}^{M-1} |\dt_x^+ u_j|^2, \quad \|u\|_\infty=\sup\limits_{j\in \mathcal{T}_M^0}|u_j|,\\
&&(u,v)=h\sum\limits_{j=1}^{M-1}u_j v_j, \quad
\langle u,v\rangle=h\sum\limits_{j=0}^{M-1}u_j v_j,\quad u, v\in X_M.
\eeas
Then it is easy to get
\be\label{innpX_M}
(-\dt_x^2 u,v)=\langle\dt_x^+ u,\dt_x^+ v\rangle,\quad
((-\dt_x^2)^{-1}u,v)=(u,(-\dt_x^2)^{-1}v),\quad u, v\in X_M.
\ee

Define the error functions $e^{\ep,k}$, $f^{\ep,k}$ as
$$e^{\ep,k}_j=E^\ep(x_j,t_k)-E_j^{\ep,k},\quad f_j^{\ep,k}=F^\ep(x_j,t_k)-F_j^{\ep,k},\quad j\in \mathcal{T}_M^0,\quad 0\le k\le \fl{T}{\tau}.$$
Then we have the following error estimates for the finite difference discretization \eqref{scheme1} with \eqref{Hd}-\eqref{E_1}.
\begin{theorem}\label{thm1}
Under the assumptions (A)-(B), there exist $h_0>0$ and $\tau_0>0$ sufficiently small and independent of $\ep$ such that, when $0<h\le h_0$, $0<\tau\le \tau_0$, the scheme \eqref{scheme1} with \eqref{Hd}-\eqref{E_1} satisfies the following error estimates
\begin{align}
\|e^{\ep,k}\|+\|\dt_x^+e^{\ep,k}\|+\|f^{\ep,k}\|&\lesssim h^2+ \fl{\tau^2}{\ep},\quad 0\le k\le \fl{T}{\tau},\quad \ep\in (0,1], \label{esti1}\\
\|e^{\ep,k}\|+\|\dt_x^+e^{\ep,k}\|+\|f^{\ep,k}\|&\lesssim h^2+\tau^2+\tau\ep^{\alpha^*}+\ep.\label{esti2}
\end{align}
Thus by taking the minimum, we have the uniform $\ep$-independent error bound
\be\label{uniform_wp}
\|e^{\ep,k}\|+\|\dt_x^+e^{\ep,k}\|+\|f^{\ep,k}\|\lesssim h^2+\tau,
\quad 0\le k\le \fl{T}{\tau},\quad \ep\in (0,1].
\ee
\end{theorem}

\section{Error analysis}
To prove Theorem \ref{thm1}, we will get the error bound \eqref{esti1} by using the energy method and \eqref{esti2} via the limiting equation KG-OP \eqref{AKG1d}, which can be displayed in the
following diagram \cite{Bao,Cai2013,Bao2014,Bao2016}.
\[
\xymatrixcolsep{8pc}\xymatrix{
(E^{\ep,k}, F^{\ep,k}) \ar[r]^{\quad O(h^2+\tau^2+
\tau\ep^{\alpha^*}+\ep^{1+\alpha^*})} \ar[rd]^{ %\qquad\quad (H^1,L^2) \text{ error}
}_{O(h^2+\tau^2/\ep)}&
(\wt{E}^\ep,0)\ar[d]^{O(\ep)} \\
&(E^\ep,F^\ep)}
\]
To simplify notations, for a function $V(x,t)$, and a grid function $V^k\in X_M$ ($k\ge 0$), we denote for $k\ge 1$
\[\Lparen V\Rparen(x,t_k)=\fl{V(x,t_{k+1})+V(x,t_{k-1})}{2},\quad x\in \overline{\Omega};\quad \Lbrack V\Rbrack_j^k=\fl{V_j^{k+1}+V_j^{k-1}}{2}, \quad j\in \mathcal{T}_M^0.\]

To bound the numerical solution, following the idea in \cite{Akrivis,Bao, Cai2013, Bao2012, Thomee}, we truncate the nonlinearity to a global Lipschitz function with compact support in $d$-dimensions, then the error can be achieved if the numerical solution is close to
the bounded exact solution. Choose a smooth function $\rho(s)\in C^\infty(\mathbb{R})$ such that
$$\rho(s)=
\left\{
\begin{aligned}
&1, \quad &|s|\le 1,\\
&\in[0,1],\quad &|s|\le 2,\\
&0,\quad &|s|\ge 2, \\
\end{aligned}
\right.
$$
and set
$$M_0=\max\left\{\sup\limits_{\ep\in(0,1]}\|E^\ep\|_{L^\infty(\Omega_T)},
\sup\limits_{\ep\in(0,1]}\|\widetilde{E}^\ep\|_{L^\infty(\Omega_T)}\right\}, $$
where $ \Omega_T=\Omega\times [0,T]$, which is well defined by assumption (A).
For $s\ge 0$, $y_1, y_2\in \mathbb{R}$, define
\be\label{rhod}
\rho_B(s)=s^2\rho(s/B),\quad B=M_0+1,
\ee
and
$$g(y_1,y_2)=\fl{1}{2}\int_0^1\rho_B'(sy_1+(1-s)y_2)ds
=\fl{\rho_B(y_1)-\rho_B(y_2)}{2(y_1-y_2)}.
$$
Then $\rho_B(s)$ is global Lipschitz and there exists $C_B>0$, such that
\be\label{rhoB}
|\rho_B(s_1)-\rho_B(s_2)|\le C_B|s_1-s_2|,\quad \forall s_1, s_2 \ge 0.
\ee
Set $\hat{E}^{\ep,0}=E^{\ep,0}$, $\hat{F}^{\ep,0}=F^{\ep,0}$,
$\hat{E}^{\ep,1}=E^{\ep,1}$, $\hat{F}^{\ep,1}=F^{\ep,1}$, and define $\hat{E}^{\ep,k}$, $\hat{F}^{\ep,k}\in X_M$ for $k\ge 1$ as following
\be\label{scheme}
\begin{aligned}
&\dt_t^2 \hat{E}_j^{\ep,k}=(\dt_x^2-1-H_j^{\ep,k})
\Lbrack\hat{E}^{\ep}\Rbrack_j^k+\big(\rho_B(\hat{E}_j^{\ep,k})
-\hat{F}^{\ep,k}_j\big)g(\hat{E}^{\ep,k+1}_j,\hat{E}_j^{\ep,k-1}),\\
&\dt_t^2\hat{F}_j^{\ep,k}=\fl{1}{2\ep^2}\dt_x^2(\hat{F}_j^{\ep,k+1}+\hat{F}_j^{\ep,k-1})+
\dt_t^2\rho_B(\hat{E}_j^{\ep,k}).
\end{aligned}
\ee
Here $(\hat{E}^{\ep,k},\hat{F}^{\ep,k})$ can be viewed as another approximation of
$(E^\ep(x_j,t_k),F^\ep(x_j,t_k))$. Applying standard fixed point arguments (refer to \cite{Cai2013}), we can get that \eqref{scheme} is uniquely solvable for sufficiently small $\tau$.

Define the error function $\hat{e}^{\ep,k}$, $\hat{f}^{\ep,k}\in X_M$ as
$$\hat{e}_j^{\ep,k}=E^\ep(x_j,t_k)-\hat{E}^{\ep,k}_j,\quad
\hat{f}_j^{\ep,k}=F^{\ep}(x_j,t_k)-\hat{F}^{\ep,k}_j, \quad j\in \mathcal{T}_M^0,\quad k\ge 0.$$
Regarding the error bounds on $(\hat{e}^{\ep,k},\hat{f}^{\ep,k})$, we have the following estimates.
\begin{theorem}\label{thm2}
Under the assumption (A), there exists $\tau_1>0$
sufficiently small, when $0<\tau\le \tau_1$, the
scheme \eqref{scheme} satisfies the following error estimates
$$\|\hat{e}^{\ep,k}\|+\|\dt_x^+\hat{e}^{\ep,k}\|+\|\hat{f}^{\ep,k}\|\lesssim h^2+ \fl{\tau^2}{\ep}, \quad 0\le k\le \fl{T}{\tau},\quad 0<\ep\le 1.$$
\end{theorem}

In order to prove it, we introduce the local truncation error $\hat{\xi}_j^{\ep,k}$, $\hat{\eta}_j^{\ep,k}\in X_M$ as
\be\label{local}
\begin{split}
\hat{\xi}_j^{\ep,k}&=\dt_t^2 E^\ep(x_j,t_k)-(\dt_x^2-1-H_j^{\ep,k})\Lparen E^\ep\Rparen(x_j,t_k)\\
&\quad-\left[\rho_B(E^\ep(x_j,t_k))- F^\ep(x_j,t_k)\right]
g\left(E^\ep(x_j,t_{k+1}),E^\ep(x_j,t_{k-1})\right)\\
&=\dt_t^2 E^\ep(x_j,t_k)-\left[\dt_x^2-1+|E^\ep(x_j,t_{k})|^2-H_j^{\ep,k}-
F^\ep(x_j,t_k)\right]\Lparen E^\ep\Rparen(x_j,t_k),\\
\hat{\eta}_j^{\ep,k}&=\ep^2\dt_t^2F^\ep(x_j,t_k)-\dt_x^2\Lparen F^\ep\Rparen(x_j,t_k)-
\ep^2\dt_t^2\rho_B(E^\ep(x_j,t_k))\\
&=\ep^2\dt_t^2F^\ep(x_j,t_k)-\dt_x^2\Lparen F^\ep\Rparen(x_j,t_k)-
\ep^2\dt_t^2|E^\ep(x_j,t_k)|^2,\quad j\in\mathcal{T}_M,\quad k\ge 1.
\end{split}
\ee
For the local truncation error, we have the following error bounds.
\begin{lemma}\label{local_e}
Under the assumption (A), we have for $j\in\mathcal{T}_M$
\[
|\hat{\xi}_j^{\ep,k}|\lesssim h^2+\fl{\tau^2}{\ep},\quad
|\hat{\eta}_j^{\ep,k}|\lesssim h^2+\tau^2,\quad 1\le k\le \fl T \tau-1;\quad
|\dt_t^c\hat{\eta}_j^{\ep,k}|\lesssim  h^2 +\fl{\tau^2}{\ep},\quad 2\le k\le \fl T \tau-2.
\]
\end{lemma}
\emph{Proof.} By \eqref{KGZ1d} and Taylor expansion, we have
\begin{align*}
&\dt_t^2 E^\ep(x_j,t_k)=\sum\limits_{m=\pm 1}\int_{0}^{1}(1-s) \p_{tt}E^\ep(x_j,t_k+ms\tau)ds\\
&=\int_{-1}^1(1-|s|)\left(\p_{xx}E^\ep-E^\ep+(E^\ep)^3
-E^\ep F^\ep-E^\ep G^\ep\right)(x_j,t_k+s\tau)ds\\
&=\p_{xx}E^\ep(x_j,t_k)-E^\ep(x_j,t_k)+E^\ep(x_j,t_k)^3- E^\ep(x_j,t_k)F^\ep(x_j,t_k)\\
&\quad+\fl{\tau^2}{6}\int_{-1}^1(1-|s|)^3\p_{tt}\left(\p_{xx}E^\ep-E^\ep+(E^\ep)^3-E^\ep F^\ep\right)(x_j,t_k+s\tau) ds\\
&\quad-\int_{-1}^1(1-|s|)E^\ep(x_j,t_k+s\tau)G^\ep\left(x_j,t_k+s\tau\right)ds.
\end{align*}
Similarly, by Taylor expansion, one can easily get that
\begin{align*}
&\left[\dt_x^2-1+|E^\ep(x_j,t_{k})|^2-H_j^{\ep,k}-F^\ep(x_j,t_k)\right]\Lparen E^\ep\Rparen(x_j,t_k)\\
&=\p_{xx}E^\ep(x_j,t_k)+\left[|E^\ep(x_j,t_{k})|^2-1-H_j^{\ep,k}-F^\ep(x_j,t_k)\right]
E^\ep(x_j,t_k)\\
&\quad+\fl{h^2}{6}\int_{-1}^1(1-|s|)^3
\Lparen \p_x^4E^\ep\Rparen(x_j+sh,t_k)ds
+\fl{\tau^2}{2}\int_{-1}^1(1-|s|)\p_x^2\p_t^2E^\ep(x_j,t_k+s\tau)ds\\
&\quad+\fl{\tau^2}{2}\left[|E^\ep(x_j,t_{k})|^2-1-H_j^{\ep,k}-F^\ep(x_j,t_k)\right]
\int_{-1}^1(1-|s|)\p_{tt}E^\ep(x_j,t_k+s\tau)ds.
\end{align*}
Note that by \eqref{Hd}, we have
\begin{align*}
&\int_{-1}^1(1-|s|)E^\ep(x_j,t_k+s\tau) G^\ep\left(x_j,t_k+s\tau\right)ds-E^\ep(x_j,t_k)H_j^{\ep,k}\\\
&=\tau E^\ep_t(x_j,t_k)\int_{0}^1s(1-s) \left[G^\ep\left(x_j,t_k+s\tau\right)-
G^\ep\left(x_j,t_k-s\tau\right)\right]ds+A_1\\
&=\tau^2E^\ep_t(x_j,t_k)\int_0^1s(1-s)\int_{-s}^s
\p_tG^\ep\left(x_j,t_k+\tht \tau\right)d\tht ds+A_1,
\end{align*}
where
$$A_1=\tau^2\int_{-1}^1(1-|s|)G^\ep\left(x_j,t_k+s\tau\right)\int_0^s (s-\tht) \p_{tt}E^\ep(x_j,t_k+\tht\tau)
d\tht ds.$$
Accordingly, by the assumption (A) and \eqref{Gp}, we deduce that
\begin{align*}
|\hat{\xi}_j^{\ep,k}|&\lesssim h^2\|\p_x^4E^\ep\|_{L^\infty}
+\tau^2 \Bigl[
\|\p_x^2\p_t^2E^\ep\|_{L^\infty}+\|\p_t^2E^\ep\|_{L^\infty}(1+\|G^\ep\|_{L^\infty}+\|F^\ep\|_{L^\infty}+
\|E^\ep\|^2_{L^\infty})\\
&\quad+\|\p_tE^\ep\|_{L^\infty}(\|\p_tG^\ep\|_{L^\infty}+\|\p_tF^\ep\|_{L^\infty})
+\|E^\ep\|_{L^\infty}(\|\p_tE^\ep\|_{L^\infty}^2
+\|\p_t^2F^\ep\|_{L^\infty})\Bigr]\\
&\lesssim h^2+\fl{\tau^2}{\ep},\quad j\in\mathcal{T}_M,\quad 1\le k\le \fl T \tau-1.
\end{align*}

Similar expansion gives
\begin{align*}
\hat{\eta}_j^{\ep,k}&=\fl{\ep^2\tau^2}{6}\int_{-1}^1(1-|s|)^3
\left[\p_t^4F^\ep(x_j,t_k+s\tau)-\p_t^4|E^\ep|^2(x_j,t_k+s\tau)\right]ds\\
&\quad-\fl{\tau^2}{2}\int_{-1}^1(1-|s|)\p_x^2\p_t^2F^\ep(x_j,t_k+s\tau)ds
-\fl{h^2}{6}\int_{-1}^1(1-|s|)^3
\Lparen\p_x^4F^\ep\Rparen(x_j+sh,t_k)ds,
\end{align*}
which implies
\begin{align*}
|\hat{\eta}_j^{\ep,k}|&\lesssim h^2\|\p_x^4F^\ep\|_{L^\infty}
+\tau^2\left(\|\p_x^2\p_t^2F^\ep\|_{L^\infty}+\ep^2\|\p_t^4F^\ep\|_{L^\infty
}+\ep^2\|\p_t^4|E^\ep|^2\|_{L^\infty}\right)\\
&\lesssim h^2+\tau^2,\quad j\in\mathcal{T}_M,\quad 1\le k\le\fl T\tau-1.
\end{align*}
Applying $\dt_t^c$ to $\hat{\eta}_j^{\ep,k}$  for $2\le k\le \fl{T}{\tau}-1$, one can deduce that
\begin{align*}
|\dt_t^c\hat{\eta}_j^{\ep,k}|&\lesssim h^2\|\p_x^4\p_tF^\ep\|_{L^\infty}
+\tau^2(\|\p_x^2\p_t^3F^\ep\|_{L^\infty}+\ep^2\|\p_t^5F^\ep\|_{L^\infty}+
\ep^2\|\p_t^5|E^\ep|^2\|_{L^\infty})\\
&\lesssim  h^2+\fl{\tau^2}{\ep},\quad j\in\mathcal{T}_M,\quad 2\le k\le\fl T\tau-2.
\end{align*}
Thus the proof is completed.
\hfill $\square$

For the initial step, we have the following estimates.
\begin{lemma} \label{initial-e}
Under the assumption (A), the first step errors of the discretization \eqref{E_1} satisfy
\[
\hat{e}_j^{\ep,0}=\hat{f}_j^{\ep,0}=0,\quad |\hat{e}_j^{\ep,1}|+|\dt_x^+\hat{e}_j^{\ep,1}|+
|\hat{f}_j^{\ep,1}|\lesssim \fl{\tau^3}{\ep},\quad
|\dt_t^+\hat{e}_j^{\ep,0}|+|\dt_t^+\hat{f}_j^{\ep,0}|\lesssim \fl{\tau^2}{\ep}.
\]
\end{lemma}
\emph{Proof.}
By the definition of $\hat{E}^{\ep,1}_j$, one can derive that
\[
|\hat{e}_j^{\ep,1}|=\fl{\tau^3}{2} \left|\int_0^1(1-s)^2\p_t^3E^\ep(x_j,s\tau)ds\right|\lesssim\tau^3\|\p_t^3E^\ep\|_{L^\infty}
  \lesssim \fl{\tau^3}{\ep},
\]
which implies that $|\dt_t^+\hat{e}_j^{\ep,0}|\lesssim \fl{\tau^2}{\ep}$.
Similarly, $|\dt_x^+\hat{e}_j^{\ep,1}|\lesssim \tau^3 \|\p_x\p_t^3 E^\ep\|_{L^\infty}\lesssim \fl{\tau^3}{\ep}$.
It follows from \eqref{E_1} and assumption (A) that
$$|\hat{f}^{\ep,1}_j|=\fl{\tau^3}{2}|\int_0^1(1-s)^2\p_t^3F^\ep(x_j,s\tau)ds|\lesssim
\tau^3\|\p_t^3F^\ep\|_{L^\infty}
  \lesssim \fl{\tau^3}{\ep}.$$
Recalling that $\hat{f}_j^{\ep,0}=0$, we can get that
$|\dt_t^+\hat{f}_j^{\ep,0}|\lesssim \fl{\tau^2}{\ep}$, which completes the proof.
\hfill $\square$

Subtracting \eqref{scheme} from \eqref{local}, we have the error equations for $j \in \mathcal{T}_M, 1\le k< \fl{T}{\tau}$,
\be\label{n_eq}
\begin{split}
&\dt_t^2\hat{e}^{\ep,k}_j=(\dt_x^2-1-H_j^{\ep,k})
\fl{\hat{e}_j^{\ep,k+1}+\hat{e}_j^{\ep,k-1}}{2}+r_j^k+\hat{\xi}_j^{\ep,k},\\
&\dt_t^2\hat{f}_j^{\ep,k}=\fl{1}{2\ep^2}\dt_x^2(\hat{f}_j^{\ep,k+1}+
\hat{f}_j^{\ep,k-1})+\dt_t^2p_j^k+\hat{\eta}_j^{\ep,k},
\end{split}
\ee
where
\be\label{pr}
\begin{split}
r_j^k=&\left(|E^\ep|^2
-F^\ep\right)\Lparen E^\ep\Rparen(x_j,t_k)-\left(\rho_B(\hat{E}_j^{\ep,k})
-\hat{F}_j^{\ep,k}\right) g(\hat{E}_j^{\ep,k+1},\hat{E}_j^{\ep,k-1}),\\
p_j^k=&|E^\ep(x_j,t_k)|^2-\rho_B(\hat{E}_j^{\ep,k}).
\end{split}
\ee
By the property of $\rho_B$ (cf. \eqref{rhoB}), one can easily get that
\be\label{p_k}
|p_j^k|=|\rho_B(E^\ep(x_j,t_k))-\rho_B(\hat{E}_j^{\ep,k})|\le C_B |\hat{e}_j^{\ep,k}|, \quad j \in \mathcal{T}_M,\quad 0\le k\le \fl T\tau.
\ee
By the definition of $g(\cdot,\cdot)$, and noting that $\Lparen E^\ep\Rparen(x_j,t_k)=g\big(E^\ep(x_j,t_{k+1}), E^\ep(x_j,t_{k-1})\big)$, it is known from \cite{Cai2015} that
for $j\in \mathcal{T}_M$, $1\le k\le\fl T\tau-1$,
\be\label{nonlbd}
\left|g(\hat{E}^{\ep,k+1}_j,\hat{E}_j^{\ep,k-1})\right|\lesssim 1, \quad  \left|\Lparen E^\ep\Rparen(x_j,t_k)-g(\hat{E}_j^{\ep,k+1},\hat{E}_j^{\ep,k-1})\right|
\lesssim \sum_{l=k\pm1} |\hat{e}_j^{\ep,l}|.
\ee
\emph{Proof of Theorem \ref{thm2}.}
Multiplying both sides of the first equation of \eqref{n_eq} by $2\tau\dt_t^c\hat{e}_j^{\ep,k}$, summing together for $j\in \mathcal{T}_M$, we obtain for $1\le k\le \fl T\tau-1$,
\be\label{e_eq2}
\begin{aligned}
&\|\dt_t^+\hat{e}^{\ep,k}\|^2-\|\dt_t^+\hat{e}^{\ep,k-1}\|^2+
\fl{1}{2}(\|\dt_x^+\hat{e}^{\ep,k+1}\|^2-\|\dt_x^+\hat{e}^{\ep,k-1}\|^2+
\|\hat{e}^{\ep,k+1}\|^2-\|\hat{e}^{\ep,k-1}\|^2)\\
&=(-H^{\ep,k}\Lbrack\hat{e}^\ep\Rbrack^k+r^k+\hat{\xi}^{\ep,k},
\hat{e}^{\ep,k+1}-\hat{e}^{\ep,k-1}).
\end{aligned}
\ee
For analyzing the second equation of \eqref{n_eq}, we introduce $\hat{u}^{\ep,k+1/2}\in X_M$ by
$$-\dt_x^2\hat{u}_j^{\ep,k+1/2}=\dt_t^+(\hat{f}_j^{\ep,k}-p_j^k).$$
Multiplying both sides of the second equation of \eqref{n_eq} by $\tau \ep^2(\hat{u}_j^{\ep,k+1/2}+\hat{u}_j^{\ep,k-1/2})$,
summing together for $j\in \mathcal{T}_M$, we have
\be\label{n_eq2}
\begin{aligned}
&\ep^2(\|\dt_x^+\hat{u}^{\ep,k+1/2}\|^2-\|\dt_x^+\hat{u}^{\ep,k-1/2}\|^2)
+\fl{1}{2}(\|\hat{f}^{\ep,k+1}\|^2-\|\hat{f}^{\ep,k-1}\|^2)\\
&=(\Lbrack\hat{f}^\ep\Rbrack^k, p^{k+1}-p^{k-1})
+\tau(\hat{\eta}^{\ep,k},\hat{u}^{\ep,k+1/2}+\hat{u}^{\ep,k-1/2}),\quad 1\le k\le \fl T\tau-1.
\end{aligned}
\ee
Introduce a discrete `energy' by
\be\label{energy_d}
\begin{split}
&\mathcal {A}^k=\|\dt_t^+\hat{e}^{\ep,k}\|^2+\fl{1}{2}
(\|\hat{e}^{\ep,k}\|^2+
\|\hat{e}^{\ep,k+1}\|^2+\|\dt_x^+\hat{e}^{\ep,k}\|^2+\|\dt_x^+
\hat{e}^{\ep,k+1}\|^2)\\
&\qquad\, +\ep^2\|\dt_x^+\hat{u}^{\ep,k+1/2}\|^2
+\fl{1}{2}(\|\hat{f}^{\ep,k+1}\|^2+\|\hat{f}^{\ep,k}\|^2),\quad 0\le k\le \fl T\tau-1.
\end{split}
\ee
Combining \eqref{e_eq2} and \eqref{n_eq2}, we get for $1\le k\le \fl T\tau-1$
\be\label{sp}
\begin{split}
&\mathcal{A}^k-\mathcal{A}^{k-1}=(-H^{\ep,k}\Lbrack\hat{e}^\ep\Rbrack^k+r^k+\hat{\xi}^{\ep,k},
\hat{e}^{\ep,k+1}-\hat{e}^{\ep,k-1})\\
&\qquad\qquad\qquad+(\Lbrack\hat{f}^\ep\Rbrack^k, p^{k+1}-p^{k-1})+\tau(\hat{\eta}^{\ep,k},\hat{u}^{\ep,k+1/2}+\hat{u}^{\ep,k-1/2}).
\end{split}
\ee

Now we estimate the terms in \eqref{sp} respectively.
By the definition of $r^k$, it can be induced that
\be\label{r_d}
\begin{split}
&r^k_j=\left(|E^\ep(x_j,t_{k})|^2
-F^\ep(x_j,t_k)\right)\left(\Lparen E^\ep\Rparen(x_j,t_k)-
g(\hat{E}^{\ep,k+1}_j,\hat{E}_j^{\ep,k-1})\right)\\
&\qquad+g(\hat{E}^{\ep,k+1}_j,\hat{E}_j^{\ep,k-1})\left(p_j^{k}
-\hat{f}_j^{\ep,k}\right).
\end{split}
\ee
In view of assumption (A), \eqref{p_k} and \eqref{nonlbd}, we derive that
\be\label{r_bound}
|r_j^k|\lesssim |\hat{e}_j^{\ep,k+1}|+|\hat{e}_j^{\ep,k}|+|\hat{e}_j^{\ep,k-1}|+
|\hat{f}^{\ep,k}_j|. \ee
This implies that
\bea
&&(-H^{\ep,k}\Lbrack\hat{e}^\ep\Rbrack^k+r^k+\hat{\xi}^{\ep,k},
\hat{e}^{\ep,k+1}-\hat{e}^{\ep,k-1})\nn\\
&&=\tau(-H^{\ep,k}\Lbrack\hat{e}^\ep\Rbrack^k+r^k+\hat{\xi}^{\ep,k},
\dt_t^+\hat{e}^{\ep,k}+\dt_t^+\hat{e}^{\ep,k-1})\nn\\
&&\lesssim \tau(1+\|H^{\ep,k}\|_\infty)(\|r^k\|^2+\|\hat{\xi}^{\ep,k}\|^2+
\sum\limits_{l=k\pm1}\|\hat{e}^{\ep,l}\|^2+\sum\limits_{l=k-1}^k
\|\dt_t^+\hat{e}^{\ep,l}\|^2)\nn\\
&&\lesssim \tau(\|\hat{\xi}^{\ep,k}\|^2+\mathcal{A}^k+\mathcal{A}^{k-1}).\label{es1_r}
\eea
It can be easily get from \eqref{nonlbd} and assumption (A) that
\begin{align*}
&p^{k+1}-p^{k-1}=E^\ep(x_j,t_{k+1})^2-E^\ep(x_j,t_{k-1})^2-4\tau g(\hat{E}_j^{\ep,k+1},
\hat{E}_j^{\ep,k-1})\dt_t^c\hat{E}_j^{\ep,k}\\
&=2(\Lparen E^\ep\Rparen(x_j,t_k)-g(\hat{E}_j^{\ep,k+1},
\hat{E}_j^{\ep,k-1}))(E^\ep(x_j,t_{k+1})-E^\ep(x_j,t_{k-1}))\\
&\quad+2g(\hat{E}_j^{\ep,k+1},\hat{E}_j^{\ep,k-1})(\hat{e}_j^{\ep,k+1}-\hat{e}_j^{\ep,k-1})\\
&\lesssim\tau \|E_t^\ep\|_{L^\infty}(|\hat{e}_j^{\ep,k+1}|+|\hat{e}_j^{\ep,k-1}|)
+\tau(|\dt_t^+\hat{e}_j^{\ep,k}|+|\dt_t^+\hat{e}_j^{\ep,k-1}|),
\end{align*}
which yields for $1\le k\le\fl T\tau-1$,
\be
(\Lbrack\hat{f}^\ep\Rbrack^k,2\tau\dt_t^c p^{k})\lesssim \tau\big[\sum\limits_{l=k-1}^k\|\dt_t^+\hat{e}^{\ep,l}\|^2+
\sum\limits_{l=k\pm1}(\|\hat{e}^{\ep,l}\|^2
+\|\hat{f}^{\ep,l}\|^2)\big]\lesssim \tau(\mathcal{A}^k+\mathcal{A}^{k-1}).\label{est-p}
\ee
Hence it can be concluded from \eqref{sp}, \eqref{es1_r} and \eqref{est-p} that
\be\label{sp2}
\mathcal{A}^k-\mathcal{A}^{k-1}-\tau(\hat{\eta}^{\ep,k},\hat{u}^{\ep,k+1/2}+
\hat{u}^{\ep,k-1/2})\lesssim\tau\left(\|\hat{\xi}^{\ep,k}\|^2
+\mathcal{A}^k+\mathcal{A}^{k-1}\right).
\ee
Applying \eqref{innpX_M}, Sobolev inequality and Cauchy inequality, we obtain
\bea\label{u_e}
&&-\fl{\mathcal{A}^k}{4}+\tau\sum\limits_{l=1}^k\left(\hat{\eta}^{\ep,l},\hat{u}^{\ep,l+1/2}+
\hat{u}^{\ep,l-1/2}\right)\nonumber\\
&&=-\fl{\mathcal{A}^k}{4}+\sum\limits_{l=1}^k\left((-\dt_x^2)^{-1}\hat{\eta}^{\ep,l},
\hat{f}^{\ep,l+1}-p^{l+1}-(\hat{f}^{\ep,l-1}-p^{l-1})\right)\nn\\
&&=-\fl{\mathcal{A}^k}{4}-2\tau\sum\limits_{l=2}^{k-1}\left(\dt_t^c(-\dt_x^2)^{-1}\hat{\eta}^{\ep,l},
\hat{f}^{\ep,l}-p^{l}\right)\nn\\
&&\quad +\sum\limits_{l=k}^{k+1}\left((-\dt_x^2)^{-1}\hat{\eta}^{\ep,l-1},\hat{f}^{\ep,l}-p^{l}\right)
-\sum\limits_{l=0}^1\left((-\dt_x^2)^{-1}\hat{\eta}^{\ep,l+1},\hat{f}^{\ep,l}-p^{l}\right)\nn\\
&&\lesssim \mathcal{A}^0+
\tau\sum\limits_{l=2}^{k-1} (\|\dt_t^c\hat{\eta}^{\ep,l}\|^2+\mathcal{A}^l)+
\sum\limits_{l=1}^2\|\hat{\eta}^{\ep,l}\|^2
+\sum\limits_{l=k-1}^{k}\|\hat{\eta}^{\ep,l}\|^2.
\eea
Summing the equation \eqref{sp2} together for $k=1,2,\cdots,m\le \fl{T}{\tau}-1$, applying \eqref{u_e}, we obtain that
\be\label{sp3}
\mathcal{A}^m\lesssim \mathcal{A}^0+\tau\sum\limits_{l=1}^{m}\mathcal{A}^l
+\sum\limits_{l=1}^2\|\hat{\eta}^{\ep,l}\|^2
+\sum\limits_{l=m-1}^{m}\|\hat{\eta}^{\ep,l}\|^2+\tau\sum\limits_{l=1}^m\|\hat{\xi}^{\ep,l}\|^2+
\tau\sum\limits_{l=2}^{m-1}\|\dt_t^c\hat{\eta}^{\ep,l}\|^2.
\ee
By Lemma \ref{initial-e} and the discrete Sobolev inequality, we deduce that
\be\label{u_e1}
\ep\|\dt_x^+\hat{u}^{\ep,1/2}\|\lesssim \ep\|\dt_t^+(\hat{f}^{\ep,0}- p^0)\|\lesssim \ep\|\dt_t^+\hat{f}^{\ep,0}\|+\ep\|\dt_t^+ \hat{e}^{\ep,0}\|
\lesssim \tau^2,
\ee
which together with Lemma \ref{initial-e} yields that
$$\mathcal{A}^0\lesssim \tau^4/\ep^2.$$
Applying Lemma \ref{local_e} and \eqref{sp3}, it can be concluded that there exists $\tau_1>0$ such that when $\tau\le \tau_1$, we have
$$\mathcal{A}^m\lesssim \left(h^2+\fl{\tau^2}{\ep}\right)^2+\tau\sum\limits_{i=1}^{m-1}\mathcal{A}^i.$$
Applying discrete Gronwall inequality, for sufficiently small $\tau$, we can conclude that
$$\mathcal{A}^m\lesssim \left(h^2+\fl{\tau^2}{\ep}\right)^2,\quad 0\le m\le \fl{T}{\tau}-1,$$
which completes the proof of Theorem \ref{thm2} by recalling \eqref{energy_d}.
\hfill $\square$

\begin{theorem}\label{thm3}
Under the assumptions (A)-(B), there exists $\tau_2>0$
sufficiently small, when $0<\tau\le \tau_2$ and $0<h\le\fl 12$, the
scheme \eqref{scheme} satisfies the following error estimates
$$\|\hat{e}^{\ep,k}\|+\|\dt_x^+\hat{e}^{\ep,k}\|+\|\hat{f}^{\ep,k}\|\lesssim h^2+\tau^2+\tau\ep^{\alpha^*}+\ep^{1+\alpha^*}, \quad 0\le k\le \fl{T}{\tau}.$$
\end{theorem}

Define another error function
$$\wt{e}^{\ep,k}_j=\wt{E}^\ep(x_j,t_k)-\hat{E}_j^{\ep,k},\quad
\wt{f}_j^{\ep,k}=-\hat{F}_j^{\ep,k},\quad j\in \mathcal{T}_M,\quad
0\le k\le \fl T\tau,$$
where $\wt{E}^\ep(x,t)$ is the solution of the KG-OP \eqref{AKG1d}. The local truncation error $\wt{\xi}^{\ep,k}$, $\wt{\eta}^{\ep,k}\in X_M$ is defined as
\be\label{local2}
\begin{split}
\wt{\xi}_j^{\ep,k}&=\dt_t^2 \wt{E}^\ep(x_j,t_k)-(\dt_x^2-1-H_j^{\ep,k})\Lparen\wt{E}^\ep\Rparen(x_j,t_k)\\
&\quad-\rho_B(\wt{E}^\ep(x_j,t_k))
g\big(\wt{E}^\ep(x_j,t_{k+1}),\wt{E}^\ep(x_j,t_{k-1})\big)\\
&=\dt_t^2 \wt{E}^\ep(x_j,t_k)-\left(\dt_x^2-1+|\wt{E}^\ep(x_j,t_k)|^2-H_j^{\ep,k}\right)
\Lparen\wt{E}^\ep\Rparen(x_j,t_k),\\
\wt{\eta}_j^{\ep,k}&=-\ep^2\dt_t^2\rho_B(\wt{E}^\ep(x_j,t_k))=-\ep^2\dt_t^2|\wt{E}^\ep(x_j,t_k)|^2.
\end{split}
\ee
\begin{lemma}\label{local2e}
Under the assumption (A), we can obtain the following error bounds
$$\|\wt{\xi}^{\ep,k}\|\lesssim h^2+\tau^2+\tau\ep^{\alpha^*},\quad
\|\wt{\eta}^{\ep,k}\|\lesssim \ep^2,\quad \|\dt_t^c\wt{\eta}^{\ep,k}\|\lesssim \ep^{1+\alpha^*}.$$
\end{lemma}
\emph{Proof.} Similar to the proof of Lemma \ref{local_e}, we can get that
\begin{align*} \wt{\xi}^{\ep,k}_j=&-\fl{h^2}{6}\int_{-1}^1(1-|s|)^3 \Lparen
\p_x^4\wt{E}^\ep \Rparen(x_j+sh,t_k)ds\\
&+\fl{\tau^2}{6}\int_{-1}^1(1-|s|)^2\p_{tt}\left[\p_{xx}\wt{E}^\ep-\wt{E}^\ep+|\wt{E}^\ep|^3\right]
(x_j,t_k+s \tau)ds\\
&-\fl{\tau^2}{2}\int_{-1}^1(1-|s|)\p_x^2\p_t^2\wt{E}^\ep(x_j,t_k+s\tau)ds-A_2\\
&-\fl{\tau^2}{2}\left(|\wt{E}^\ep(x_j,t_{k})|^2-1-H_j^{\ep,k}\right)
\int_{-1}^1(1-|s|)\p_{tt}\wt{E}^\ep(x_j,t_k+s\tau)ds,
\end{align*}
where
\begin{align*}
A_2&=\int_{-1}^1(1-|s|) \wt{E}^\ep(x_j,t_k+s\tau) G^\ep\left(x_j,t_k+s\tau\right)ds-\wt{E}^\ep(x_j,t_k)H_j^{\ep,k}\\
&=\tau\int_{-1}^1(1-|s|)G^\ep\left(x_j,t_k+s\tau\right)\int_0^s
\p_t\wt{E}^\ep(x_j,t_k+\tht\tau)d\tht ds\\
&\lesssim \tau \|G^\ep\|_{L^\infty}\|\p_t\wt{E}^\ep\|_{L^\infty}
\lesssim \tau\ep^{\alpha^*}.
\end{align*}
Hence we can conclude from assumption (A) that
$$\|\wt{\xi}^{\ep,k}\|\lesssim h^2+\tau^2+\tau\ep^{\alpha^*}.$$
Note that by assumption (A), it is easy to get that
$$\p_t^3|\wt{E}^\ep|^2= 6\p_t\wt{E}^\ep\p_{tt}\wt{E}^\ep+2\wt{E}^\ep\p_t^3\wt{E}^\ep
\lesssim \ep^{\alpha^*-1},$$
which indicates that
$$\|\wt{\eta}^{\ep,k}\|\lesssim \ep^2,\quad
\|\dt_t^c \wt{\eta}^{\ep,k}\|\lesssim \ep^{1+\alpha^*},$$
the proof is completed.
\hfill $\square$

Analogous to Lemma \ref{initial-e}, we have the error bounds for $\wt{e}^{\ep,k}$, $\wt{f}^{\ep,k}$ at the first step.
\begin{lemma} \label{initial-e2}
Under the assumptions (A) and (B), the first step errors of the discretization \eqref{E_1} satisfy
\[
\wt{e}_j^{\ep,0}=\wt{f}_j^{\ep,0}=0,\quad  |\wt{e}_j^{\ep,1}|+|\dt_x^+\wt{e}_j^{\ep,1}|\lesssim \tau^3+\tau^2\ep^{\alpha^*},\quad
|\wt{f}_j^{\ep,1}|+|\dt_t^+\wt{e}_j^{\ep,0}|\lesssim \tau^2+\tau \ep^{\alpha^*},\quad
|\dt_t^+\wt{f}_j^{\ep,0}|\lesssim \tau.
\]
\end{lemma}
\emph{Proof.} It follows from \eqref{AKG} and \eqref{KGZ1d} that
$\p_{tt}E^\ep(x_j,0)=\p_{tt}\wt{E}^\ep(x_j,0)$. By \eqref{E_1} and assumption (B), one gets that
\begin{align*}
\wt{e}_j^{\ep,1}&=\fl{\tau^3}{2}\int_0^1(1-s)^2\p_t^3\wt{E}^\ep(x_j,s\tau)ds\\
&=\fl{\tau^3}{2}\int_0^1(1-s)^2\p_t\left(\p_{xx}\wt{E}^\ep-\wt{E}^\ep+
|\wt{E}^\ep|^3-\wt{E}^\ep G^\ep\right)(x_j,s\tau)ds\\
&=\fl{\tau^3}{2}\int_0^1(1-s)^2\p_t\left(\p_{xx}\wt{E}^\ep-\wt{E}^\ep+
|\wt{E}^\ep|^3\right)(x_j,s\tau)ds\\
&\quad+\fl{\tau^2}{2}E_0(x_j)\ep^\alpha \og_0(x_j)-\tau^2\int_0^1(1-s)\wt{E}^\ep(x_j,s\tau)G^\ep(x_j,s\tau)ds\\
&\lesssim \tau^3+\tau^2\ep^{\alpha^*}.
\end{align*}
Thus this gives that $|\dt_t^+\wt{e}_j^{\ep,0}|\lesssim \tau^2+\tau\ep^{\alpha^*}$. Similar arguments can deduce that
$|\dt_x^+\wt{e}_j^{\ep,1}|\lesssim \tau^3+\tau^2\ep^{\alpha^*}$.
By the definition, we have
$$|\wt{f}^{\ep,1}_j|=|F_j^{\ep,1}|\lesssim \tau^2 |\p_{tt}F^\ep(x_j,0)|\lesssim \tau^2.$$
The remaining conclusions are direct.
\hfill $\square$

\emph{Proof of Theorem \ref{thm3}.}
Subtracting \eqref{scheme} from \eqref{local2}, one has the error equations
\be\label{2e_eq}
\begin{split}
&\dt_t^2\wt{e}^{\ep,k}_j=(\dt_x^2-1-H_j^{\ep,k})\Lbrack\wt{e}^\ep\Rbrack_j^k
+\wt{r}_j^k+\wt{\xi}_j^{\ep,k},\\
&\ep^2\dt_t^2\wt{f}_j^{\ep,k}=\dt_x^2\Lbrack\wt{f}^\ep\Rbrack_j^k
+\ep^2\dt_t^2\wt{p}_j^k+\wt{\eta}_j^{\ep,k},\quad j\in \mathcal{T}_M,\quad 1\le k\le \fl T\tau-1,
\end{split}
\ee
where
\begin{align*}
\wt{r}_j^k&=|\wt{E}^\ep(x_j,t_k)|^2\Lparen\wt{E}^\ep\Rparen(x_j,t_k)
-\left(\rho_B(\hat{E}_j^{\ep,k})-\Lbrack\hat{F}^\ep\Rbrack^k_j\right) g(\hat{E}_j^{\ep,k+1},\hat{E}_j^{\ep,k-1}),\\
\wt{p}_j^k&=|\wt{E}^\ep(x_j,t_k)|^2-\rho_B(\hat{E}_j^{\ep,k}),
\end{align*}
Suppose $\wt{u}^{\ep,k+\fl{1}{2}}\in X_M$ is the solution to the equation
$$-\dt_x^2\wt{u}_j^{\ep,k+\fl{1}{2}}=\dt_t^+(\wt{f}_j^{\ep,k}-\wt{p}_j^k),\quad j\in\mathcal{T}_M,\quad 0\le k\le \fl T\tau-1.$$
Denote
\be\label{energy_d2}
\begin{aligned}
\wt{\mathcal {A}}^k=&\|\dt_t^+\wt{e}^{\ep,k}\|^2+\fl{1}{2}\left(\|\wt{e}^{\ep,k}\|^2+
\|\wt{e}^{\ep,k+1}\|^2+\|\dt_x^+\wt{e}^{\ep,k}\|^2+\|\dt_x^+\wt{e}^{\ep,k+1}\|^2\right)\nonumber\\
&+\ep^2\|\dt_x^+\wt{u}^{\ep,k+1/2}\|^2
+\fl{1}{2}\left(\|\wt{f}^{\ep,k}\|^2+\|\wt{f}^{\ep,k+1}\|^2\right).
\end{aligned}
\ee
Applying the same approach as in the former part, there exists $\tau_2>0$ sufficiently small independent of $\ep$ such that
\[
\wt{\mathcal{A}}^k\lesssim \wt{\mathcal{A}}^0+\tau\sum\limits_{l=1}^{k-1}\wt{\mathcal{A}}^l
+\sum\limits_{l=1}^2\|\wt{\eta}^{\ep,l}\|^2
+\sum\limits_{l=k-1}^{k}\|\wt{\eta}^{\ep,l}\|^2+\tau\sum\limits_{l=1}^k
\|\wt{\xi}^{\ep,l}\|^2+
\tau\sum\limits_{l=2}^{k-1}\|\dt_t\wt{\eta}^{\ep,l}\|^2.
\]
By Lemma \ref{initial-e2} and the discrete Sobolev inequality, we deduce that
$$\ep\|\dt_x^+\wt{u}^{\ep,1/2}\|\lesssim \ep\|\dt_t^+\wt{f}^{\ep,0}\|+\ep
\|\dt_t^+\wt{e}^{\ep,0}\|\lesssim \ep\tau.$$
which together with Lemma \ref{initial-e2} yields that
$$\wt{\mathcal{A}}^0\lesssim (\tau^2+\tau\ep^{\alpha^*}+\ep\tau)^2.$$
Applying Lemma \ref{local2e}, it can be concluded that when $0<\tau\le\tau_2$,
$$\wt{\mathcal {A}}^k\lesssim (h^2+\tau^2+\ep\tau+\tau\ep^{\alpha^*}
+\ep^{1+\alpha^*})^2+\tau\sum\limits_{i=1}^{k-1}\wt{\mathcal{A}}^i.$$
It follows from discrete Gronwall inequality that
$$\wt{\mathcal{A}}^k\lesssim (h^2+\tau^2+\ep\tau+\tau\ep^{\alpha^*}+
\ep^{1+\alpha^*})^2,$$
implying that
$$\|\wt{e}^{\ep,k}\|+\|\dt_x^+\wt{e}^{\ep,k}\|+\|\wt{f}^{\ep,k}\|\lesssim
h^2+\tau^2+\ep\tau+\tau\ep^{\alpha^*}+\ep^{1+\alpha^*}.$$
Using the assumption (B) and the triangle inequality, we obtain that
\begin{align*}
\|\hat{e}^{\ep,k}\|+\|\dt_x^+\hat{e}^{\ep,k}\|&\lesssim\|\wt{e}^{\ep,k}\|
+\|\dt_x^+\wt{e}^{\ep,k}\|+\|E^\ep(\cdot,t_k)-\wt{E}^\ep(\cdot,t_k)\|_{H^1}\\
&\lesssim h^2+\tau^2+\tau\ep^{\alpha^*}+\ep,\\
\|\hat{f}^{\ep,k}\|&\lesssim \|\wt{f}^{\ep,k}\|+\|F^\ep(\cdot,t_k)\|_{L^2}
\lesssim h^2+\tau^2+\tau\ep^{\alpha^*}+\ep,
\end{align*}
which completes the proof of Theorem \ref{thm3}.
\hfill $\square$

\bigskip

\emph{Proof of theorem \ref{thm1}.}
Now we have proved the two types of estimates \eqref{esti1} and \eqref{esti2} for ($\hat{E}^{\ep,k},\hat{F}^{\ep,k}$), which is the
solution of the modified finite difference discretization \eqref{scheme} with \eqref{Hd} and \eqref{E_1}. Hence we can get the
uniform error bounds for $(\hat{E}^{\ep,k},\hat{F}^{\ep,k})$:
$$\|\hat{e}^{\ep,k}\|+\|\dt_x^+\hat{e}^{\ep,k}\|+\|\hat{f}^{\ep,k}\|
\lesssim h^2+\tau,$$
which together with the inverse inequality \cite{Thomee} yields
\[\|\hat{E}^{\ep,k}\|_\infty-\|E^\ep(\cdot,t_k)\|_{L^\infty}\le \|\hat{e}^{\ep,k}\|_\infty\lesssim \|\dt_x^+\hat{e}^{\ep,k}\|\lesssim h^2+\tau, \quad 0\le k\le \fl T\tau.\]
Thus there exists $h_0>0$ and $\tau_3>0$ sufficiently small such that when $0<h\le h_0$ and $0<\tau\le \tau_3$,
\[\|\hat{E}^{\ep,k}\|_\infty\le 1+\|E^\ep(\cdot,t_k)\|_\infty\le 1+M_0,\quad 0\le k\le \fl T\tau.\]
Set $\tau_0=\min\{\tau_1,\tau_2,\tau_3\}$, when $0<h\le h_0$, $0<\tau\le \tau_0$, \eqref{scheme} collapses to \eqref{scheme1}, i.e. ($\hat{E}^{\ep,k},\hat{F}^{\ep,k}$) are identical to ($E^{\ep,k},F^{\ep,k}$), which completes the proof.
\hfill $\square$

\begin{remark}
The error bounds in Theorem \ref{thm1} are still valid in higher dimensions, e.g.
$d=2,3$. The key point is the discrete Sobolev inequality in higher dimensions as \cite{Cai2013,Thomee}
$$\|\psi_h\|_\infty\le \fl{1}{C_d(h)}\|\psi_h\|_{H^1},\quad \mathrm{with}\quad
 C_d(h)\sim
 \left\{
\begin{aligned}
&\fl{1}{|\mathrm{ln}\, h|},\quad &d=2,\\
&h^{1/2},\quad &d=3. \\
\end{aligned}
\right.$$
where $\psi_h$ is a mesh function over $\Omega$ with homogeneous Dirichlet boundary condition.
Thus by requiring an additional condition on the time step $\tau$
\[\tau=o(C_d(h)),\]
the same error bounds can be obtained.

\end{remark}

\section{Numerical results}
In this section, we present numerical results for the KGZ \eqref{KGZ1d} by
the finite difference discretization \eqref{scheme1} with \eqref{Hd}-\eqref{E_1}. In our experiment, the initial condition is set as
\beas
&&E_0(x)=e^{-x^2}\sin x,\quad E_1(x)=\mathrm{sech}(x^2/2)\cos x,\\
&& \og_0(x)=\mathrm{sech} (x^2)\cos(3x),\quad \og_1(x)=\mathrm{sech}( x^2)\sin(4x),
\eeas
and the parameters $\alpha$ and $\beta$ are chosen as

Case I. $\alpha=1$ and $\beta=0$;

Case II. $\alpha=0$ and $\beta=-1$.
\medskip

\begin{table}[h!]
\def\temptablewidth{1\textwidth}
\vspace{-12pt}
\caption{Spatial errors at time $t=1$ for Case II, i.e. $\alpha=0$, $\beta=-1$.}\label{ill-h}
{\rule{\temptablewidth}{1pt}}
\begin{tabular*}{\temptablewidth}{@{\extracolsep{\fill}}ccccccc}
$e^\ep(1)$&$h_0=0.2$&$h_0/2$&$h_0/2^2$&$h_0/2^3$&$h_0/2^4$&$h_0/2^5$\\\hline
$\ep=1$&	1.57E-2&	4.05E-3&	1.02E-3&	2.56E-4&	6.39E-5&	1.60E-5\\
rate&-&	1.95&	1.99&	2.00&	2.00&	2.00\\\hline
$\ep=1/2$&	1.35E-2&	3.48E-3&	8.76E-4&	2.19E-4&	5.49E-5&	1.37E-5\\
rate&-&	1.95&	1.99&	2.00&	2.00&	2.00\\\hline
$\ep=1/2^{2}$&	1.30E-2&	3.35E-3&	8.44E-4&	2.11E-4&	5.29E-5&	1.32E-5\\
rate&-&	1.95&	1.99&	2.00&	2.00&	2.00\\\hline
$\ep=1/2^{3}$&	1.32E-2&	3.42E-3&	8.60E-4&	2.15E-4&	5.39E-5&	1.35E-5\\
rate&-&	1.95&	1.99&	2.00&	2.00&	2.00\\\hline
$\ep=1/2^{4}$&	1.33E-2&	3.43E-3&	8.65E-4&	2.17E-4&	5.42E-5&	1.36E-5\\
rate&-&	1.95&	1.99&	2.00&	2.00&	2.00\\\hline
$\ep=1/2^{5}$&	1.33E-2&	3.44E-3&	8.66E-4&	2.17E-4&	5.43E-5&	1.36E-5\\
rate&-&	1.95&	1.99&	2.00&	2.00&	2.00\\\hline
$\ep=1/2^{6}$&	1.33E-2&	3.44E-3&	8.66E-4&	2.17E-4&	5.42E-5&	1.36E-5\\
rate&-&	1.95&	1.99&	2.00&	2.00&	2.00\\\hline
$\ep=1/2^{7}$&	1.33E-2&	3.44E-3&	8.65E-4&	2.17E-4&	5.42E-5&	1.36E-5\\
rate&-&	1.95&	1.99&	2.00&	2.00&	2.00\\\hline
$\ep=1/2^{8}$&	1.33E-2&	3.43E-3&	8.65E-4&	2.17E-4&	5.42E-5&	1.36E-5\\
rate&-&	1.95&	1.99&	2.00&	2.00&	2.00\\
 \toprule
  \bottomrule
$n^\ep(1)$&$h_0=0.2$&$h_0/2$&$h_0/2^2$&$h_0/2^3$&$h_0/2^4$&$h_0/2^5$\\\hline
$\ep=1$&	1.91E-2&	4.79E-3&	1.20E-3&	2.99E-4&	7.49E-5&	1.87E-5\\
rate&-&	2.00&	2.00&	2.00&	2.00&	2.00\\\hline
$\ep=1/2$&		1.61E-2&	3.98E-3&	9.92E-4&	2.48E-4&	6.20E-5&	1.55E-5\\
rate&-&	2.02&	2.00&	2.00&	2.00&	2.00\\\hline
$\ep=1/2^{2}$&	6.59E-3&	1.67E-3&	4.18E-4&	1.05E-4&	2.62E-5&	6.56E-6\\
rate&-&	1.98&	2.00&	2.00&	2.00&	1.99\\\hline
$\ep=1/2^{3}$&	5.30E-3&	1.35E-3&	3.39E-4&	8.49E-5&	2.13E-5&	5.33E-6\\
rate&-&	1.97&	1.99&	2.00&	2.00&	2.00\\\hline
$\ep=1/2^{4}$&	5.12E-3&	1.30E-3&	3.28E-4&	8.20E-5&	2.05E-5&	5.15E-6\\
rate&-&	1.97&	1.99&	2.00&	2.00&	1.99\\\hline
$\ep=1/2^{5}$&	5.06E-3&	1.29E-3&	3.23E-4&	8.10E-5&	2.03E-5&	5.09E-6\\
rate&-&	1.97&	1.99&	2.00&	2.00&	1.99\\\hline
$\ep=1/2^{6}$&	5.03E-3&	1.28E-3&	3.21E-4&	8.05E-5&	2.02E-5&	5.08E-6\\
rate&-&	1.97&	1.99&	2.00&	2.00&	1.99\\\hline
$\ep=1/2^{7}$&	5.01E-3&	1.28E-3&	3.21E-4&	8.02E-5&	2.01E-5&	5.07E-6\\
rate&-&	1.97&	1.99&	2.00&	2.00&	1.99\\\hline
$\ep=1/2^{8}$&	5.01E-3&	1.27E-3&	3.20E-4&	8.01E-5&	2.01E-5&	5.07E-6\\
rate&-&	1.97&	1.99&	2.00&	1.99&	1.99\\
\end{tabular*}
{\rule{\temptablewidth}{1pt}}
\end{table}

\begin{table}[h!]
\def\temptablewidth{1\textwidth}
\setlength{\tabcolsep}{3pt}
\vspace{-12pt}
\caption{Temporal errors at time $t=1$ for Case I, i.e. $\alpha=1$, $\beta=0$.}\label{well-t}
{\rule{\temptablewidth}{1pt}}
\begin{tabular*}{\temptablewidth}{@{\extracolsep{\fill}}ccccccccc}
$e^\ep(1)$&$\tau_0=0.05$&$\tau_0/2$&$\tau_0/2^2$&$\tau_0/2^3$&$\tau_0/2^4$&$\tau_0/2^5$&
$\tau_0/2^6$&$\tau_0/2^7$\\\hline
$\ep=1$&	6.34E-3&	1.64E-3&	4.16E-4&	1.05E-4&	2.62E-5&	6.59E-6&	1.67E-6& 4.20E-7\\
rate&-&	1.95&	1.98&	1.99&	1.99&	1.99&	1.98& 1.99\\\hline
$\ep=1/2$&	5.19E-3&	1.35E-3&	3.42E-4&	8.61E-5&	2.16E-5&	5.43E-6&	1.37E-6& 3.45E-7\\
rate&-&	1.95&	1.98&	1.99&	1.99&	1.99&	1.98&1.99\\\hline
$\ep=1/2^{2}$&	5.17E-3&	1.34E-3&	3.41E-4&	8.58E-5&	2.15E-5&	5.41E-6&	1.37E-6& 3.45E-7\\
rate&-&	1.95&	1.98&	1.99&	1.99&	1.99&	1.98& 1.99\\\hline
$\ep=1/2^{3}$&	5.17E-3&	1.34E-3&	3.41E-4&	8.58E-5&	2.15E-5&	5.41E-6&	1.37E-6& 3.45E-7\\
rate&-&	1.95&	1.98&	1.99&	1.99&	1.99&	1.98& 1.99\\\hline
$\ep=1/2^{4}$&	5.16E-3&	1.34E-3&	3.40E-4&	8.56E-5&	2.15E-5&	5.40E-6&	1.37E-6&  3.45e-7\\
rate&-&	1.95&	1.98&	1.99&	1.99&	1.99&	1.98& 1.99\\\hline
$\ep=1/2^{5}$&	5.14E-3&	1.34E-3&	3.39E-4&	8.55E-5&	2.15E-5&	5.39E-6&	1.36E-6&	3.44E-7\\
rate&-&	1.94&	1.98&	1.99&	1.99&	1.99&	1.98&	1.98\\\hline
$\ep=1/2^{6}$&	5.13E-3&	1.33E-3&	3.39E-4&	8.54E-5&	2.15E-5&	5.39E-6&	1.36E-6&	3.45E-7\\
rate&-&	1.95&	1.97&	1.99&	1.99&	1.99&	1.98&	1.98\\\hline
$\ep=1/2^{7}$&	5.00E-3&	1.30E-3&	3.29E-4&	8.32E-5&	2.09E-5&	5.26E-6&	1.33E-6&	3.36E-7\\
rate&-&	1.95&	1.98&	1.98&	1.99&	1.99&	1.98&	1.98\\\hline
$\ep=1/2^{8}$&	5.01E-3&	1.30E-3&	3.29E-4&	8.29E-5&	2.09E-5&	5.25E-6&	1.33E-6&	3.36E-7\\
rate&-&	1.95&	1.98&	1.99&	1.99&	1.99&	1.98&	1.98\\
\toprule
\bottomrule
$n^\ep(1)$&$\tau_0=0.05$&$\tau_0/2$&$\tau_0/2^2$&$\tau_0/2^3$&$\tau_0/2^4$&$\tau_0/2^5$&
$\tau_0/2^6$&$\tau_0/2^7$\\\hline
$\ep=1$&	6.59E-3&	1.70E-3&	4.29E-4&	1.08E-4&	2.70E-5&	6.78E-6&	1.71E-6& 4.28E-7\\
rate&-&	1.96&	1.98&	1.99&	2.00&	2.00&	1.99& 2.00\\\hline
$\ep=1/2$&	1.87E-2&	4.85E-3&	1.23E-3&	3.08E-4&	7.71E-5&	1.93E-5&	4.85E-6& 1.21E-6\\
rate&-&	1.95&	1.98&	1.99&	2.00&	2.00&	1.99& 2.00\\\hline
$\ep=1/2^{2}$&	1.50E-2&	4.63E-3&	1.24E-3&	3.15E-4&	7.90E-5&	1.98E-5&	4.96E-6& 1.24E-6\\
rate&-&	1.70&	1.90&	1.98&	1.99&	2.00&	2.00& 2.00\\\hline
$\ep=1/2^{3}$&	8.65E-3&	3.50E-3&	1.44E-3&	4.36E-4&	1.13E-4&	2.83E-5&	7.09E-6& 1.77E-6\\
rate&-&	1.31&	1.28&	1.73&	1.95&	1.99&	2.00& 2.00\\\hline
$\ep=1/2^{4}$&	5.55E-3&	1.96E-3&	9.14E-4&	4.87E-4&	1.86E-4&	5.02E-5&	1.27E-5& 3.18E-6\\
rate&-&	1.50&	1.10&	0.91&	1.39&	1.89&	1.99& 2.00\\\hline
$\ep=1/2^{5}$&	4.97E-3&	1.45E-3&	5.53E-4&	2.66E-4&	1.56E-4&	7.76E-5&	2.38E-5&	6.09E-6\\
rate&-&	1.78&	1.39&	1.06&	0.77&	1.01&	1.70&	1.97\\\hline
$\ep=1/2^{6}$&	4.70E-3&	1.30E-3&	4.38E-4&	1.81E-4&	8.58E-5&	4.86E-5&	2.86E-5&	1.12E-5\\
rate&-&	1.85&	1.57&	1.27&	1.08&	0.82&	0.76&	1.36\\\hline
$\ep=1/2^{7}$&	4.19E-3&	1.20E-3&	3.57E-4&	1.41E-4&	6.30E-5&	2.96E-5&	1.55E-5&	9.45E-6\\
rate&-&	1.81&	1.75&	1.34&	1.16&	1.09&	0.93&	0.72\\\hline
$\ep=1/2^{8}$&	3.96E-3&	1.12E-3&	3.28E-4&	1.12E-4&	4.84E-5&	2.26E-5&	1.07E-5&	5.30E-6\\
rate&-&	1.83&	1.77&	1.55&	1.21&	1.10&	1.08&	1.01\\
\end{tabular*}
{\rule{\temptablewidth}{1pt}}
\end{table}

\begin{table}[h!]
\def\temptablewidth{1\textwidth}
\setlength{\tabcolsep}{3pt}
\vspace{-12pt}
\caption{Temporal errors at time $t=1$ for Case II, i.e. $\alpha=0$, $\beta=-1$.}\label{ill-t}
{\rule{\temptablewidth}{1pt}}
\begin{tabular*}{\temptablewidth}{@{\extracolsep{\fill}}ccccccccc}
$e^\ep(1)$&$\tau_0=0.05$&$\tau_0/2$&$\tau_0/2^2$&$\tau_0/2^3$&$\tau_0/2^4$&$\tau_0/2^5$&
$\tau_0/2^6$&$\tau_0/2^7$\\\hline
$\ep=1$&6.34E-3&	1.64E-3&	4.16E-4&	1.05E-4&	2.62E-5&	6.59E-6&	1.67E-6&  4.20E-7\\
rate&-&1.95&	1.98&	1.99&	1.99&	1.99&	1.98& 1.99\\\hline
$\ep=1/2$&	5.48E-3&	1.42E-3&	3.62E-4&	9.11E-5&	2.29E-5&	5.74E-6&	1.45E-6& 3.65E-7\\
rate&-&	1.94&	1.98&	1.99&	1.99&	1.99&	1.98& 1.99\\\hline
$\ep=1/2^{2}$&5.16E-3&	1.34E-3&	3.41E-4&	8.58E-5&	2.16E-5&	5.41E-6&	1.37E-6&  3.42E-7\\
rate&-&	1.95&	1.98&	1.99&	1.99&	1.99&	1.98& 1.99\\\hline
$\ep=1/2^{3}$&5.20E-3&	1.35E-3&	3.44E-4&	8.66E-5&	2.18E-5&	5.46E-6&	1.38E-6&  3.47E-7\\
rate&-&1.94&	1.98&	1.99&	1.99&	1.99&	1.98& 1.99\\\hline
$\ep=1/2^{4}$&5.47E-3&	1.42E-3&	3.59E-4&	9.04E-5&	2.27E-5&	5.71E-6&	1.44E-6& 3.63E-7\\
rate&-&1.95&	1.98&	1.99&	1.99&	1.99&	1.98& 1.99\\\hline
$\ep=1/2^{5}$&5.77E-3&	1.63E-3&	4.04E-4&	1.00E-4&	2.51E-5&	6.30E-6&	1.59E-6&4.01E-7\\
rate&-&1.83&	2.01&	2.01&	2.00&	1.99&	1.99& 1.99\\\hline
$\ep=1/2^{6}$&	5.52E-3&	1.98E-3&	5.67E-4&	1.32E-4&	3.15E-5&	7.81E-6&	1.96E-6& 4.93E-7\\
rate&-&	1.48&	1.80&	2.11&	2.06&	2.01&	1.99&1.99\\\hline
$\ep=1/2^{7}$&	5.45E-3&	1.92E-3&	8.40E-4&	2.40E-4&	5.20E-5&	1.18E-5&	2.88E-6&	7.26E-7\\
rate&-&1.51&	1.19&	1.81&	2.21&	2.13&	2.04&	1.99\\\hline
$\ep=1/2^{8}$&	5.45E-3&	1.90E-3&	8.42E-4&	4.03E-4&	1.13E-4&	2.35E-5&	5.17E-6&	1.24E-6\\
rate&-&	1.52&	1.17&	1.06&	1.83&	2.27&	2.19&	2.06\\
\toprule
\bottomrule
$n^\ep(1)$&$\tau_0=0.05$&$\tau_0/2$&$\tau_0/2^2$&$\tau_0/2^3$&$\tau_0/2^4$&$\tau_0/2^5$&
$\tau_0/2^6$&$\tau_0/2^7$\\\hline
$\ep=1$&	6.59E-3&	1.70E-3&	4.29E-4&	1.08E-4&	2.70E-5&	6.78E-6&	1.71E-6& 4.28E-7\\
rate&-&	1.96&	1.98&	1.99&	2.00&	2.00&	1.99&2.00\\\hline
$\ep=1/2$&1.33E-2&	3.45E-3&	8.73E-4&	2.19E-4&	5.49E-5&	1.38E-5&	3.46E-6& 8.65E-7\\
rate&-&1.95&	1.98&	1.99&	2.00&	2.00&	1.99& 2.00\\\hline
$\ep=1/2^{2}$&	9.61E-3&	2.97E-3&	7.98E-4&	2.03E-4&	5.09E-5&	1.28E-5&	3.20E-6& 8.00E-7\\
rate&-&1.70&	1.89&	1.98&	1.99&	2.00&	2.00& 2.00\\\hline
$\ep=1/2^{3}$&5.55E-3&	2.30E-3&	9.03E-4&	2.77E-4&	7.21E-5&	1.81E-5&	4.54E-6& 1.14E-6\\
rate&-&1.27&	1.35&	1.70&	1.94&	1.99&	2.00& 2.00\\\hline
$\ep=1/2^{4}$&3.92E-3&	1.35E-3&	6.30E-4&	3.13E-4&	1.17E-4&	3.23E-5&	8.18E-6& 2.05E-6\\
rate&-&1.54&	1.10&	1.01&	1.43&	1.85&	1.98& 2.00\\\hline
$\ep=1/2^{5}$&3.91E-3&	1.19E-3&	4.18E-4&	1.89E-4&	1.07E-4&	4.80E-5&	1.53E-5&	3.96E-6\\\
rate&-&	1.71&	1.51&	1.14&	0.82&	1.16&	1.65&	1.95\\\hline
$\ep=1/2^{6}$&	3.53E-3&	1.37E-3&	4.38E-4&	1.49E-4&	6.27E-5&	3.45E-5&	1.84E-5&	7.03E-6\\
rate&-&	1.37&	1.64&	1.55&	1.25&	0.86&	0.91&	1.39\\\hline
$\ep=1/2^{7}$&		3.31E-3&	1.27E-3&	5.49E-4&	1.78E-4&	5.74E-5&	2.26E-5&	1.11E-5&	6.51E-6\\
rate&-&	1.38&	1.21&	1.63&	1.63&	1.35&	1.02&	0.77\\\hline
$\ep=1/2^{8}$&		3.18E-3&	1.20E-3&	5.47E-4&	2.52E-4&	7.81E-5&	2.34E-5&	8.60E-6&	3.80E-6\\
rate&-&	1.41&	1.13&	1.12&	1.69&	1.74&	1.45&	1.18\\
\end{tabular*}
{\rule{\temptablewidth}{1pt}}
\end{table}

In practical computation, the truncated domain is set as $\Omega_\varepsilon=\left[-30-\frac{1}{\varepsilon},
30+\frac{1}{\varepsilon}\right]$, which is large enough such that the homogeneous Dirichlet boundary condition does not introduce significant errors. Similar to the truncation for the Zakharov system, the bounded computational domain $\Omega_\varepsilon$
has to be chosen as $\varepsilon$-dependent due to that the rapid outgoing waves are at wave speed $O\left(\frac{1}{\varepsilon}\right)$ and the homogeneous Dirichlet boundary condition is taken at the boundary. The computational $\ep$-dependent domain can be fixed as
$\ep$-independent if one applies absorbing boundary condition (ABC) \cite{Eng} or transport boundary condition (TBC) \cite{Feng,Giv},
or perfected matched layer (PML) \cite{Bere} for the wave-type equations in \eqref{KGZ1d} and \eqref{wave} during the truncation (refer to \cite{Bao2016}).

To quantify the numerical errors, we introduce the error functions as follows
$$e^\ep(t_k):=\fl{\|e^{\ep,k}\|+\|\dt_x^+e^{\ep,k}\|}{\|E^\ep(\cdot,t_k)\|_{H^1}},\quad n^\ep(t_k):=\fl{\|n^{\ep,k}\|}{\|N^\ep(\cdot,t_k)\|_{L^2}},$$
where $e^{\ep,k}=E^\ep(\cdot,t_k)-E^{\ep,k}$, $n^{\ep,k}=N^\ep(\cdot,t_k)-N^{\ep,k}$.
The ``exact" solution is obtained by the EWI-SP method \cite{Cai2013} with very small mesh size $h=1/64$ and time step $\tau=10^{-6}$.
The errors are displayed at $t=1$. For spatial error analysis, we set a time step $\tau=10^{-5}$, such that the temporal error can be neglected; for temporal error analysis, the mesh size $h$ is set as $h=2.5\times 10^{-4}$ such that the spatial error can be ignored.

Table \ref{ill-h} depicts the spatial errors for  Case II initial data, which clearly
demonstrates that our numerical method is uniformly second order accurate in $h$ for all $\ep \in (0,1]$. The result for Case I initial data is similar, which is omitted here for brevity.

Tables \ref{well-t} and \ref{ill-t} present the temporal errors for
Cases I and II, respectively, from which we can conclude that the method is uniformly convergent in time for both initial data. Specifically, Table \ref{well-t} shows the method is uniformly second order accurate for $E^\ep$, while for $N^\ep$, it is second
order in time when $\tau\lesssim \ep$ or $\ep\lesssim \tau^2$ (cf. upper and lower
triangle parts, respectively). There is a resonance regime when $\tau\sim \ep$ where the
convergence rate degenerates to the first order, which agrees with the analysis
\eqref{esti1}-\eqref{esti2}. For $\alpha=0$, $\beta=-1$, the upper and lower triangle parts of Table \ref{ill-t} suggest that the method is second and first order in time when $\tau\lesssim \ep$ and $\ep\lesssim \tau$, respectively. Moreover, the upper triangle parts of Tables \ref{well-t} and \ref{ill-t} show the order of the errors at $O(\tau^2/\ep)$ for $n^\ep$ (cf. each column), which confirms our error analysis in Section 3.

\section{Conclusion}
We presented a uniformly accurate finite difference method and carried out its rigorous error bounds for the Klein-Gordon Zakharov (KGZ) system in $d$ ($d=1, 2, 3$) dimensions, which involves a dimensionless
parameter $\ep \in (0,1]$. When $0<\ep\ll1$, i.e. subsonic limit regime,
the solution of KGZ propagates highly oscillatory waves in time and/or
rapid outgoing waves in space. Our method was designed by reformulating KGZ into an asymptotic consistent formulation followed by adopting an integral approximation
for the oscillating term. By applying the energy method and
the limiting equation, two independent error bounds were obtained, which depend explicitly on the parameter $\ep$, mesh size $h$ and time step $\tau$. Thus it can be established that the method is uniformly convergent for $\ep\in(0,1]$ with quadratic and linear convergence in space and time, respectively. The error bounds is confirmed by the numerical
results, which also suggest that our estimates are sharp.

% ===========================================
%                                                                           REFERENCES
% =============================================================================

\end{document}